\documentclass[a4paper,11pt]{article}

\usepackage{amsmath,amssymb,url,graphicx, amsthm}
\usepackage[T2A]{fontenc}
\usepackage[utf8]{inputenc}

\usepackage{mathtext}
\usepackage[T2C]{fontenc}
\usepackage{mathrsfs}
\usepackage{bbm}
\usepackage{mathtools,setspace}
\usepackage{fullpage}
\usepackage[russian,english]{babel}
\usepackage{authblk} 
\usepackage{verbatim}

\title{Markov birth-and-death dynamics of populations}

\author{Viktor Bezborodov}
\affil{Fakultät für Mathematik, Universität Bielefeld, 33501 Bielefeld, Germany}

\newcommand{\R}{\mathbb{R}}
\newcommand{\N}{\mathbb{N}}
\newcommand{\Z}{\mathbb{Z}}
\newcommand{\G}{\CYRG(\R^{d})}
\newcommand{\Go}{\CYRG _0(\R^{d})}

\newcommand{\X}{\mathbf{X}}

\swapnumbers
\theoremstyle{plain}
\newtheorem{thm}{Theorem}[section]
\newtheorem{lem1}[thm]{Theorem}
\newtheorem{lem3}[thm]{Proposition}
\newtheorem{lem4}[thm]{Lemma}
\newtheorem{lem5}[thm]{Lemma}
\newtheorem{lem6}[thm]{Lemma}
\newtheorem{lem7}[thm]{Lemma}
\newtheorem{3stm}[thm]{Statement}

\newtheorem{rmk1}[thm]{Proposition}
\newtheorem{prop1}[thm]{Proposition}

\theoremstyle{definition}

\newtheorem{rmk2}[thm]{Remark}

\newtheorem{rmk5}[thm]{Remark}

\newtheorem{rmk7}[thm]{Remark}
\newtheorem{rmk8}[thm]{Remark}
\newtheorem{rmk9}[thm]{Remark}
 \newtheorem{def1}[thm]{Definition}

 \newtheorem{def5}[thm]{Definition}

\theoremstyle{plain}
\newtheorem{3thm}{Theorem}[section]
\newtheorem{3lem1}[3thm]{Proposition}
\newtheorem{3lem2}[3thm]{Proposition}
\newtheorem{3thm2}[3thm]{Theorem}

\newtheorem{3thm6}[3thm]{Theorem}
\newtheorem{3cor1}[3thm]{Corollary}
\newtheorem{3prop1}[3thm]{Proposition}
\newtheorem{3thm5}[3thm]{Theorem}
\newtheorem{3thm3}[3thm]{Theorem}
\newtheorem{3Assumpt}[3thm]{Continuity assumptions}

\theoremstyle{definition}

\newtheorem{def mart pr}[3thm]{Definition}
\newtheorem{def weak solution}[3thm]{Definition}
\newtheorem{def strong solution}[3thm]{Definition}
\newtheorem{def pathwise uniqueness}[3thm]{Definition}
\newtheorem{def joint uniqueness in law}[3thm]{Definition}
\newtheorem{3crl mart pr}[3thm]{Corollary}
\newtheorem{3rmk2}[3thm]{Remark}
\newtheorem{3rmk3}[3thm]{Remark}
\newtheorem{3rmk4}[3thm]{Remark}
\newtheorem{3rmk5}[3thm]{Remark}
\newtheorem{3rmk6}[3thm]{Remark}
\newtheorem{3rmk7}[3thm]{Remark}
\newtheorem{3rmk8}[3thm]{Remark}
\newtheorem{3rmk9}[3thm]{Remark}
\newtheorem{3rmk10}[3thm]{Remark}

\theoremstyle{plain}

\newtheorem{8prop1}[thm]{Proposition}
\newtheorem{8prop2}[thm]{Proposition}

\theoremstyle{definition}
\newtheorem{8rmk10}[thm]{Remark}

\begin{document}

\maketitle


\begin{abstract}

 Spatial birth-and-death processes with a finite number of particles
 are obtained as unique solutions to
 certain stochastic equations.
 Conditions are given for existence and uniqueness of such solutions,
 as well as for continuous dependence on the initial conditions.
 The possibility of an explosion and connection with
 the heuristic generator of the process are discussed.

\end{abstract}

\tableofcontents

\section*{Introduction}
\addcontentsline{toc}{section}{Introduction}

This article deals with spatial birth-and-death processes
which may describe stochastic dynamics of spatial
population. 
Specifically, at each moment of time
the population is 
represented as a collection 
of motionless points in 
 $\R^d$. We interpret 
the points as particles, or individuals.
Existing particles 
may die and new particles may appear.
Each particle is characterized by its 
location.

 The state space of a spatial birth-and-death
 Markov process on $\R ^d$ with finite number of points
 is the space of finite configurations over $\R ^d$, 
 \[
 \CYRG _0(\R ^d)=\{ \eta \subset \R ^d : |\eta| < \infty \},
\]
where $|\eta|$ is the number of points of $\eta$.


 Denote by $\mathscr{B}(\R ^d)$ the Borel $\sigma$-algebra on $\R ^d$.
 The evolution of a spatial birth-and-death process in $\R^d$
 admits the following description. Two functions characterize the 
 development in time,
 the birth rate coefficient
 $b: \R ^d \times \Gamma _0 (\R^d) \rightarrow [0;\infty)$
 and the death rate coefficient
 $d: \R ^d \times \Gamma _0 (\R^d) \rightarrow [0;\infty)$. 
 If the system is in state $\eta \in \Go$ 
 at time $t$, then the probability 
 that a new particle appears (a ``birth'') in a bounded set $B\in \mathscr{B}(\R ^d)$
 over time interval $[t;t+ \Delta t]$ is 
 \[
   \Delta t \int\limits _{ B}b(x, \eta)dx + o(\Delta t),
 \]
 the probability that
a particle $x \in \eta$ is deleted from the configuration (a ``death'') over 
 time interval $[t;t+ \Delta t]$ is 
 \[
  d(x,\eta) \Delta t + o(\Delta t),
 \]
 and no two events happen simultaneously.  
 By an event we mean a birth or a death.
Using a slightly different terminology, we can say that
the rate at which a birth 
 occurs in $B$ is $\int_B b(x, \eta)dx$,
 the rate at which a particle $x\in \eta$
dies is $d(x,\eta)$, and no two events happen at the same time.

 Such processes, in which the birth and death rates
depend on the spatial structure of the system as opposed to
 classical $\Z _+$-valued birth-and-death processes (see e.g. \cite{KarlinMcGregor},
 \cite{Onaseminal}, \cite[Page 116]{Branch2},
 \cite[Page 109]{Branch1}, and references therein),
 were first
 studied by Preston in \cite{Preston}. 
 A heuristic description
 similar to that above
 appeared already there. Our description 
 resembles the one in \cite{GarciaKurtz}.

The (heuristic) generator of a spatial
birth-and-death process should be of the form

\begin{align} \label{intr generator}
L F (\eta) = \int\limits _{x \in \R^d} b(x, \eta) [F(\eta \cup {x}) -F(\eta)] dx + 
\sum\limits _{x \in \eta} d(x, \eta ) (F(\eta \setminus {x}) - F(\eta)),
\end{align}
for
 $F$ in an appropriate domain, where $\eta \cup {x}$
 and $\eta \setminus {x}$ are shorthands for 
 $\eta \cup \{ x \}$ and
 $\eta \setminus \{ x \}$, respectively.

  Spatial point processes have been used in statistics for
 simulation purposes, see e.g. \cite{Statsim}, \cite[chapter 11]{Simbook}
 and references therein. For application of spatial and stochastic models in biology see
 e.g. \cite{Levin}, \cite{Ecology}, and references therein.

 To construct a spatial birth-and-death process with given birth and 
 death rate coefficients,
 we consider in Section 2 stochastic equations with Poisson type noise

\begin{equation} \label{intr se}
\begin{split}
\eta _t (B) = \int\limits _{B \times (0;t] \times [0; \infty ] }
I _{ [0;b(x,\eta _{s-} )] } (u) dN_1(x,s,u) \\
 - \int\limits _{\Z \times (0;t] \times [0; \infty ) } I_{\{x_i \in \eta _{r-} \cap B \}}
I _{ [0;d(x_i,\eta _{r-} )] } (v)  d  N _2 (i,r,v)
\end{split}
\end{equation}
where
$(\eta _t)_{t \geq 0}$ is a suitable $\Go$-valued cadlag 
stochastic process, the ``solution'' of the equation, $I_A$
is the indicator function of the set $A$, 
$B \in \mathscr{B} (\R ^d) $ is a Borel set, $N_1$
is a  Poisson point processes on $ \R^d \times \R _+ \times \R _+ $
with intensity $dx \times ds \times du $,
$  N_2$ is a Poisson point process on 
$ \Z  \times \R _+ \times \R _+$ with intensity $\# \times dr \times dv$,
$\# $ is the counting measure on $ \Z^d$,
$ \eta _0$
is a (random) initial finite  configuration,
 $b,d: \R ^d \times \Gamma _0 (\R^d) \rightarrow [0;\infty)$ are functions
 that are
measurable with respect to the product $\sigma$-algebra 
$\mathscr{B} (\R) \times \mathscr{B} ( \CYRG _0 (\R)) $ and
$\{x_i \}$ is some collection of points
satisfying $ \eta _s \subset \{x_i \}$
for every moment of time $s$
(the precise definition is given in Section \ref{voranbringen}).
We require the
processes $N_1, N_2, \eta _0$ to be independent of each other. Equation \eqref{intr se}
is understood in the sense that the equality
holds a.s. for all bounded
$B \in \mathscr{B} (\R ^d) $ 
and $t \geq 0$.

Garcia and Kurtz studied in  \cite{GarciaKurtz}
equations similar to \eqref{intr se} for infinite systems. In 
the earlier work \cite{Garcia} of Garcia  another approach was used:
birth-and-death processes were obtained
as projections of Poisson point processes. A
further development of the projection method
appears in \cite{GarciaKurtz2}.
Fournier and Meleard in \cite{FournierMeleard} considered a similar 
equation for the construction of the Bolker-Pacala-Dieckmann-Law process
with finitely many particles.

Holley and Stroock \cite{HolleyStroock} constructed a
spatial birth-and-death process
as a Markov family of  unique solutions to the corresponding martingale problem.
For the most part, they consider a process
contained in a bounded volume, with bounded 
birth and death rate coefficients.  They also
 proved the corresponding result
 for the nearest neighbor model in $\R ^1$
 with an infinite number of particles.

 Kondratiev and Skorokhod \cite{KondSkor}
 constructed a
  contact process in continuum, with the infinite  number of particles.
  The contact process can be described 
  as a spatial birth-and-death process with
  \[
   b(x,\eta)  = \lambda \sum\limits _{y \in \eta} a(x-y), 
   \ \ \  d(x,\eta) \equiv 1,
  \]
where $\lambda >0$ and $0\leq a \in L ^1 (\R ^d)$.
  Under some additional assumptions, 
  they showed existence of
  the process for a broad class of 
  initial conditions. Furthermore, if 
  the value of some 
  energy functional on 
  the initial condition is finite, 
  then it stays finite at any point in time.

 In the aforementioned references as well as in the present work
  the evolution of the system in time via Markov process is described.
  An alternative approach consists in using the concept
  of statistical dynamics that substitutes the notion of a Markov
stochastic process. This approach is based on considering evolutions of measures 
and their correlation functions. For details see e.g. 
\cite{BaDdynamics}, \cite{Statdynamics}, and references therein.

  There is an enormous amount of literature concerning
  interacting particle systems on lattices and related topics
  (e.g., \cite{Liggett}, \cite{Liggett2}, 
  \cite{Scalinglimits}, \cite{IPS1}, \cite{IPS2}, \cite{SpitzerBaD}, etc.)
Penrose in \cite{Penrose} gives a general existence 
result for interacting particle systems on a lattice
with local interactions and bounded 
jump rates (see also \cite[Chapter 9]{Liggett}).
The spin space is allowed
to be non-compact, which gives the opportunity 
to incorporate spatial 
birth-and-death processes in continuum. 
Unfortunately, the assumptions
become rather restrictive when applied 
to continuous space models. More specifically, 
the birth rate coefficient should be bounded, and 
for every bounded Borel set $B$
the expression
\[
 \sum\limits _{x \in \eta \cap B} d(x,\eta)
\]
should be bounded uniformly in $\eta$, $\eta \in \G$.

Let us briefly describe  the contents of the article.

In Section 1 we introduce give some general notions, definitions
and results related to Markov processes in configuration spaces. 
We start with configuration spaces, 
which are the state spaces for birth-and-death processes, then we
introduce and discuss metrical and topological structures thereof.
Also, we present some facts and constructions from 
probability theory, such as integration
with respect to a Poisson point process, or a
sufficient condition for a functional transformation 
of a Markov chain to be a Markov chain again.

In the second section we construct a spatial
birth-and-death process $(\eta _t)_{t \geq 0}$ as a unique
solution to equation
\eqref{intr se}. We prove strong existence and 
pathwise 
uniqueness for \eqref{intr se}.
A key condition is that we require 
$b$ to grow not faster than 
linearly in the sense that

\begin{equation} \label{intr sublinear growth for b} 
\int\limits _{\R ^d} {b}(x, \eta ) dx \leq c_1|\eta| +c_2.
\end{equation}

The equation is solved pathwisely,
``from one jump to another''. 
Also, we prove 
uniqueness 
in law for equation \eqref{intr se}
and the Markov property for the unique solution.
Considering 
\eqref{intr se} with a (non-random) initial condition
$\alpha \in \Go$ and denoting corresponding solution
by $(\eta (\alpha , t))_{t \geq 0}$, we see that 
a unique solution induces a Markov family
of probability measures on the Skorokhod space 
$D_{\Go}[0;\infty)$ (which can be regarded as 
the canonical space for a solution of \eqref{intr se}).

When the birth and death rate coefficients $b$ and $d$
satisfy some continuity assumptions, the solution
is expected to have continuous dependence on 
the initial condition, at least in some proper sense.
Realization of this idea and precise formulations
are given in Section 2.1. 
The proof is based on considering a coupling of two
birth-and-death processes.

The formal relation of a unique solution to  \eqref{intr se}
and operator $L$ in \eqref{intr generator} is 
given via the martingale problem, in Section 2.2,
and via some kind of a pointwise convergence, in Section
2.5.

In Section 2.4 we formulate and prove a
theorem about coupling of two 
birth-and-death processes.
The idea to compare a spatial birth-and-death process
with some ``simpler'' process goes back
to Preston, \cite{Preston}. In \cite{FournierMeleard}
this technique was applied to the study of the probability 
of extinction.

\section{Configuration spaces and Markov processes: miscellaneous}

 In this section we list some notions and facts we use in this work.

  \subsection{Some notations and conventions }

 Sometimes we write $\infty$ and $+ \infty$ interchangeably, so that
 $f \to \infty$ and $f \to +\infty$, or 
 $a < \infty$ and $a < + \infty$ may have the same meaning. However,
 $+ \infty$ is reserved for the real line only, whereas $\infty$
 have wider range of applications, e.g. for a sequence
 $\{ x_n \}_{n \in \N} \subset \R ^d$ we
 may write $x_n \to \infty$, $n \to \infty$, which is
 equivalent to $|x_n| \to + \infty$.
 On the other hand, we do not assign
 any meaning to $x_n \to + \infty$.


 In all probabilistic constructions we work on some
 probability space $(\Omega , \mathscr{F}, P)$,
 sometimes equipped with a filtration of $\sigma$-algebras.
Elements of $\Omega $ are usually denoted as $\omega$.

 The set $A^c$ is the complement of  the set $A \subset \Omega$: 
 $A^c = \Omega \setminus A $. We write $[a;b]$, $[a;b)$ etc. for 
 the intervals of real numbers. For example, 
 $(a;b] = \{x\in \R \mid a <x \leq b \}$, $-\infty \leq a <b \leq + \infty$.
 The half line $\R _+$ includes $0$: $\R _+ = [0; \infty)$.

 \subsection{Configuration spaces}\label{Configurations}
 
 In this section we introduce notions and facts
about spaces of configurations,
in particular, topological and
 metrical structures on $\G$ as well as a
 characterization of compact sets of $\G$.
We discuss configurations over Euclidean spaces only.

\begin{def1}
For $d\in\mathbb{N}$ and a measurable set $\Lambda\subset\mathbb{R}^d$, 
the configuration space $\CYRG(\Lambda)$ is defined as 
\[
 \CYRG(\Lambda)=\{ \gamma \subset \Lambda : |\gamma \cap K| <
 + \infty \text{ for any compact } K \subset  \R^d \}.
\]
\end{def1}
We recall that $|A|$ denotes the number of elements 
of A. We also say that $\CYRG(\Lambda)$ is the space of configurations
over $\Lambda$.
Note that $\varnothing \in \CYRG(\Lambda)$.

 Let $\Z _+$ be the set $\{0,1,2,... \}$.
 We say that a Radon measure $\mu$ on $(\R ^d, \mathscr{B}(\R ^d))$
 is a \textit{counting measure} 
 on $\R^d$ if
 $\mu (A) \in \Z_+$ for all $A \in \mathscr{B} (\R^d)$. 
 When a counting measure $\nu$ satisfies
 additionally $\nu (\{x \}) \leq 1$ for all $x \in \R ^d$, we call
 it a \textit{simple counting measure}.

 As long as it does not lead to ambiguities, we identify \label{identify conf meas}
  a configuration
 with a simple counting Radon measures
 on $\R^d$: 
 as a measure, a configuration $\gamma \in \CYRG(\R^{d})$
 maps a set $B\in \mathscr{B}$ into $|\gamma \cap B|$. In 
 other words, $\gamma = \sum\limits _{x \in \gamma} \delta _x$.

 One equips $\G$ with the vague topology, i.e., the weakest topology such that 
 for all $f\in C_c (\R^d)$ (the set of continuous functions
 on $\R^d$ with compact support)
 the map 
 \begin{align*}
\G \ni \gamma \mapsto \langle \gamma , f \rangle
:= \sum\limits _{x \in \gamma} f(x) \in \R
\end{align*}
 is continuous.

 Equipped with this topology, $\G$ is a Polish space, i.e., 
there exists a metric on $\G$ compatible with the vague topology and 
 with respect to which $\G$ is a complete separable metric space, see, e.g., \cite{KondKut},
 and references therein. 
 We say that a metric is compatible with a given topology
 if the topology induced by the 
 metric coincides with the given topology.
 
 


For a bounded $B \subset \R^d$ and $\gamma \in \G$, we denote 
$\delta ( \gamma, B) = min \{ |x-y|: x,y \in \gamma \cap B, x \neq y \} $.
Let $B_r(x)$ denote the closed ball in  $\R^d$ of 
the radius $r$ centered at  $x $.

A set is said to be \emph{relatively compact} if its closure is compact.
The following theorem gives a characterization of compact sets in $\G$, cf.
\cite{KondKut}, \cite{HolleyStroock}.

\begin{lem1} \label{compactness in G}
A set $ F \subset \CYRG (\R^d)$ is relatively compact in 
the vague topology
if and only if 

\begin{equation} \label{comp condition}
 \sup\limits_{\gamma \in F} \{ \gamma (B_n(0)) + \delta ^{-1} ( \gamma , B_n(0) ) \} < \infty 
\end{equation}
holds for all $n \in \N$.
\end{lem1}

    \textbf{Proof}. Assume that \eqref{comp condition} is 
    satisfied for some $ F \subset \CYRG (\R^d)$. In metric spaces
compactness is equivalent to sequential compactness, therefore it is
sufficient to show that an arbitrary sequence 
contains a convergent  subsequence in $\CYRG (\R^{d})$. To this end, consider an
arbitrary sequence
$ \{ \gamma _n \} _{n \in \N} \subset F$. 
The supremum $ \sup\limits_{n} \gamma _n (B_1(0)) $ is finite, consequently, 
by the Banach–Alaoglu theorem
there exists a measure $ \alpha _1 \in C  (B_1(0)) ^* $ \bigg(here $C (B_1(0)) ^* $
is the dual space of $C  (B_1(0))$ \bigg) and a subsequence $\{ \gamma _n ^{(1)} \} 
\subset \{ \gamma _n  \}$ such that 
$  \gamma _n ^{(1)} | _{B_1(0)} \to \alpha _1 $ in $ C  (B_1(0)) ^* $. 
Furthermore, one may see that $ \alpha _1 \in \CYRG (B_1(0)) $ (it is particularly
important here that $ \sup\limits_{\gamma \in F} \{ \delta ^{-1} ( \gamma , B_1(0) ) \} < \infty $ ). 
Indeed, arguing by contradiction one 
may get that $\alpha _1 (A) \in \Z _+$ for all Borel sets $A$,
and Lemma \ref{caveat} below ensures that $\alpha _1$
is a simple counting measure.

    Similarly, from the sequence  $ \gamma _n ^{(1)}$ we may extract subsequence 
$\{ \gamma _n ^{(2)}\} \subset \{ \gamma _n ^{(1)}\} $ in such a way that 
$ \gamma _n ^{(2)} $ converges to some $ \alpha _2 \in  \CYRG (B_2 (0) ) $. 
Continuing in the same way, we will find a sequence of sequences
$ \{\gamma _n ^{(m)}\}$ such that $ \gamma _n ^{(m)} \to \alpha _m \in \CYRG (B_m (0) )$ and 
$\{ \gamma _n ^{(m+1)}\} \subset \{ \gamma _n ^{(m)}\}$.
Consider now the sequence $\{ \gamma _n ^{(n)}\}_{n \in \N}$. For any $m$, restrictions of its elements to 
$B_m (0)$ converge to $ \alpha _m $ in $ \CYRG (B_m (0))$, 
Therefore, $ \gamma _n ^{(n)} \to \alpha $ in $ \CYRG (\R^d)$, where $ \alpha = \bigcup\limits_{n} \alpha _n $.

Conversely, if \eqref{comp condition} is not fulfilled 
for some $n_0 \in \N$, then we can construct
a sequence $ \{ \gamma _n \} _{n \in \N} \subset F$ such that either the first
summand in \eqref{comp condition} tends to infinity:
\[
\gamma _j (B_{n_0}(0)) \to \infty, j \to \infty
\]
in which case, of course, there is
 no convergent subsequence, or the second
summand in \eqref{comp condition} tends to infinity. 
 In the latter case,
  a subsequence of the sequence
  $ \{ \gamma _n  | _{ B_{n_0}(0)} \} _{n \in \N} $  may converge 
to a 
counting measure (when all $\gamma _n$ are considered 
as measures). However, the limit measure can not be a simple counting measure.
 Thus, the sequence
  ${\{ \gamma _n \} _{n \in \N} \subset F}$
does not contain a 
 convergent subsequence in $\G$. $\Box$

We denote by $ CS (\CYRG (\R^d)) $ the space of all compact subsets of $ \CYRG (\R^d) $.

    \begin{rmk1} \label{non sigma compactness}
The topological space $ \G $ is not $\sigma$ - compact.
     
    \end{rmk1}

    \textbf{Proof}. Let $ \{ K_m \} _{n \in \N}$ be an arbitrary sequence from $ CS (\CYRG (\R^d)) $. 
We will show that $\bigcup\limits_{n} K_n \neq \CYRG (\R^d) $.
To each compact $K_m$ we may assign
a sequence $q_1^{(m)}, q _2 ^{(m)}, ...$ of positive numbers such that 

$$ \sup\limits_{\gamma \in K_m} \{ \gamma (B_n(0)) + 
\delta ^{-1} ( \gamma , B_n(0) ) \} < q_n ^{(m)}. $$

There exists a configuration whose intersection with $B_n(0)$ contains at least 
$q_n^{(n)} +1 $ points, for each $n \in \N$.  This configuration 
does not belong to any of the sets $ \{ K_m \} _{m \in \N}$, 
hence it can not belong to the union $\bigcup\limits_{m} K_m $. $\Box$

\textbf{Remark}.
 Since $\CYRG (\R^d) $ is a separable metrizable space, 
 Proposition \ref{non sigma compactness} implies that 
$\CYRG (\R^d) $ is \textit{not locally compact}.

 For another description of all compact sets in $ \CYRG (\R^d) $ 
we will use the set $ \Phi \subset C  (\R^d) $ of all positive 
continuous functions $ \phi$ satisfying the following conditions:

 1) $ \phi (x) = \phi (y) $ whenever $ |x|=|y|$, $x, y \in \R^d$,

 2) $ \lim _{ |x| \to \infty} \phi (x) = 0 $.
 
For $ \phi \in \Phi$ we denote 
$$
\Psi = \Psi _{ \phi } (x,y) := \phi (x) \phi (y) \frac{|x-y|+1}{|x-y|} I \{ x \neq y \}.
$$

\begin{lem3}

(i) For all $c > 0$ and $ \phi \in \Phi$

$$ K_c := \bigg\{ \gamma : \iint\limits_{\R^d \times \R^d} \Psi  _{ \phi } (x,y) 
\gamma (dx) \gamma (dy) \leqslant c \bigg\} \in CS ( \CYRG (\R^d) );$$

(ii) For all $K \in CS ( \CYRG (\R^d) )$ there exist $ \phi \in \Phi$ such that

$$\sup\limits_{ \gamma \in K} ~ \{ \iint\limits_{\R^d \times \R^d} 
\Psi  _{ \phi } (x,y) \gamma (dx) \gamma (dy) \} \leqslant 1. $$

\end{lem3}
Proof. (i) Denote $ \theta _n = \min\limits_{x \in B_n (0)} \phi (x) >0 $.

For $\gamma \in K _ c$ we have 
\[
c \geqslant \iint\limits_{B_n (0) \times B_n (0)} \Psi (x,y) \gamma (dx) \gamma (dy) 
\]
\[
\geqslant \iint\limits_{B_n (0) \times B_n (0)} \phi (x) \phi (y) I \{ x \neq y \} \gamma (dx) \gamma (dy) \geqslant 
\theta_{n} ^2  \gamma (B_n (0)) ( \gamma (B_n (0)) -1 )
\]
and 
$$
c \geqslant \iint\limits_{B_n (0) \times B_n (0)} \Psi (x,y) \gamma (dx) 
\gamma (dy) \geqslant \theta_{n} ^2 \frac{\delta ^{-1} ( \gamma , B_n(0) ) +1}{\delta ^{-1} ( \gamma , B_n(0) )}
\geqslant \theta_{n} ^2 \delta ^{-1} ( \gamma , B_n(0) ).
$$
Consequently, 
\[
\sup\limits_{ \gamma \in K_c} \gamma (B_n (0))  \leqslant \theta _n \sqrt {c} +1 ,
\]
and
\[
\sup\limits_{ \gamma \in K_c} \delta ^{-1} ( \gamma , B_n(0) ) \leqslant \frac {c}{\theta _n  ^2} .
\]
It 
remains to show that $K_c$ is closed, in which case
Theorem \ref{compactness in G} will imply compactness of $ K_c$.
The space $\G$ is metrizable, therefore sequential closedness will suffice. Take 
$\gamma _k \in K _c$, $\gamma _k \to \gamma $ in $\G$, $k \to \infty$. 
For $n \in \N$, let $\Psi _n \in C_c (\R^d \times \R^d)$ be 
an increasing sequence of functions such that
$\Psi _n \leqslant \Psi$, $\Psi _n (x,y) = \Psi (x,y)$ 
for $x,y \in \R^d$ satisfying $|x|,|y| \leq n$, $|x-y|\geq \frac1n$.
For such a sequence  we have $\Psi _n (x,y) \uparrow \Psi  (x,y)$ 
 for all $x,y \in \R^d$, $x \ne y $.
For each $f \in C_c (\R^d \times \R^d)$,
the 
map 
$$
\eta \mapsto \langle \eta \times \eta , f \rangle :=
\iint\limits_{\R^d \times \R^d} f (x,y) \eta (dx) \eta (dy)
$$
is continuous 
in the vague topology. Thus for all $n\in \N$, 
$\langle \gamma _k \times \gamma _k , \Psi _n \rangle \to 
\langle \gamma  \times \gamma  , \Psi _n \rangle $.
Consequently,
$\langle \gamma  \times \gamma  , \Psi _n \rangle \leq c$,
$n \in \N$, and by Fatou's Lemma
\[
\langle \gamma  \times \gamma  , \Psi  \rangle =
\iint\limits_{\R^d \times \R^d} \Psi (x,y) \gamma (dx) \gamma (dy) =
\]
\[
=\iint\limits_{\R^d \times \R^d} \liminf\limits _n \Psi_n (x,y)
\gamma (dx) \gamma (dy) \leq \liminf\limits _n
\iint\limits_{\R^d \times \R^d}  \Psi_n (x,y)
\gamma (dx) \gamma (dy) \leq c.
\]

To prove (ii), for a given compact set $K \subset \CYRG (\R^d) $
and a given function $ \phi \in \Phi$,
denote
\[
a_n (K) := \sup\limits_{\gamma \in K} \{ \gamma (B_n(0)) + \delta ^{-1} ( \gamma , B_n(0) ) \} 
\]
and
\[
 b_n ( \phi ) := \sup\limits_{|x| > n} | \phi (x) | .
\]
Theorem 
\ref{compactness in G} implies $a_n (K) < \infty$, and we can estimate

\[
\iint\limits_{\big(B_{n+1} (0) \setminus B_n (0) \big) \times \big(B_{n+1} (0) \setminus B_n (0) \big)} 
\Psi (x,y) \gamma (dx) \gamma (dy)=
\] 
\[= \iint\limits_{\big(B_{n+1} (0) \setminus B_n (0) \big) \times \big(B_{n+1} (0) \setminus B_n (0) \big)}
\phi (x) \phi (y) \frac{|x-y|+1}{|x-y|} I \{ x \neq y \} \gamma (dx) \gamma (dy) \leqslant 
\] 
\[\leqslant
\iint\limits_{\big(B_{n+1} (0) \setminus B_n (0) \big) \times \big(B_{n+1} (0) \setminus B_n (0) \big)} 
b_n ^2 (a_n +1) \gamma (dx) \gamma (dy) \leqslant
 b_n ^2 (a_n + 1) ^3.
 \]
Taking a function $ \phi \in \Phi$  such that
\[
3 b_n ^2 ( \phi ) (a_n + 1) ^3 < 
\frac{6}{\pi ^2} \frac{1}{(n+1)^2},
\]
we get
$$
\sup\limits_{ \gamma \in K} ~ 
\{ \iint\limits_{\R^d \times \R^d} \Psi (x,y) \gamma (dx) \gamma (dy) \} \leqslant 1 . \Box 
$$

 \subsubsection{The space of finite configurations }

 For $\Lambda \subset \R^d$, the space $\CYRG _0(\Lambda)$ is 
 defined as
\[
 \CYRG _0(\Lambda) := \{ \eta \subset \Lambda : |\eta| < \infty \}.
\]
  We 
  see that $\CYRG _0(\Lambda)$ is the collection of 
 all finite subsets of $\Lambda$.
 We denote the space of $n$-point configurations as 
 $\CYRG _0 ^{(n)}(\Lambda)$:
 \[
 \CYRG _0 ^{(n)}(\Lambda) := \{ \eta \in \CYRG _0(\Lambda) \mid
 |\eta| = n\},  \ \ \ n \in \N,
\]
and
$\CYRG _0 ^{(0)}(\Lambda) := \{ \varnothing \}$. 
 Sometimes we will write $\CYRG _0$ instead of $\Go$.
 Recall that 
 we occasionally write $\eta \setminus x$ instead of 
 $\eta \setminus \{ x\} $ , $\eta \cup x$
 instead of $\eta \cup \{x\}$.
 
 To define a topological structure on $\Go$, we introduce 
 the following surjections 
 (see, e.g., \cite{KondKuna} and references therein)
 
 \begin{equation}
  \begin{split}
  sym : \bigsqcup\limits _{n=0}^\infty \widetilde {(\R^d)^n} \rightarrow \Go \\
    sym((x_1,...,x_n)) = \{x_1 ,..., x_n  \},
  \end{split}
 \end{equation}
where

\begin{equation}
\widetilde {(\R^d)^n} := \{ (x_1 ,..., x_n) \in (\R^d)^n \mid x_j \in \R ^d, j=1,...,n, 
x_i \ne x_j , i \ne j \},
\end{equation}
and, 
by convention, $\widetilde {(\R^d)^0} = \{ \varnothing \}$.

 The map $sym$ produces a one-to-one correspondence between
 $\CYRG _0 ^{(n)}(\R ^d)$, $n \geq 1$, and
 the quotient space
 $\widetilde {(\R^d)^n} / \sim _n$, where $\sim _n$ is the 
 equivalence relation on $(\R^d)^n$, 
 
 $$
 {(x_1 ,..., x_n) \sim _n (y_1 ,..., y_n)}
 $$ 
 when
 there exist a permutation
 $\sigma : \{ 1,...,n\} \rightarrow \{ 1,...,n\}$
 such that 
 
 $$
 {(x_{\sigma (1)} ,..., x_{\sigma (n)}) = (y_1 ,..., y_n)}.
 $$
 
 We endow $\CYRG _0 ^{(n)}(\R ^d)$ with the 
 topology induced by this one-to-one correspondence. 
 Equivalently, a set $A \subset \CYRG _0 ^{(n)}(\R ^d)$ is 
 open iff $sym ^{-1}(A)$ is open in $\widetilde {(\R^d)^n} $.
 The space $\widetilde {(\R^d)^n} \subset (\R^d)^n$ we consider, of course,
 with the relative, or subspace, topology.
 As far as $\CYRG _0 ^{(0)}(\R ^d) = \{ \varnothing \}$ is concerned, 
 we regard it as an open set.

Having defined topological
 structures on $\CYRG _0 ^{(n)}(\R ^d)$, $n \geq 0$,
 we endow $\Go$ with the topology of disjoint union,

 \begin{equation}
  \Go  = \bigsqcup\limits _{n=0} ^\infty \CYRG _0 ^{(n)}(\R ^d).
 \end{equation}

 In this topology, a set $K \subset \Go$ is compact iff 
 $K \subset \bigsqcup\limits _{n=0} ^N \CYRG _0 ^{(n)}(\R ^d)$
 for some $N \in \N$ and for each $n \leq N$ the set
 $K \cap \CYRG _0 ^{(n)}(\R ^d)$ is compact in 
 $\CYRG _0 ^{(n)}(\R ^d)$. A set $K_n \subset \CYRG _0 ^{(n)}(\R ^d)$
 is compact iff $sym ^{-1} (K_n)$ is compact in $\widetilde {(\R^d)^n}$.
 We note that 
 in order for $K_n$ to be compact, the set $sym ^{-1}K _n$, 
 regarded as a subset of $(\R^d)^n$,
 should not have limit points on the diagonals, i.e.
 limit points from the set $(\R^d)^n \setminus \widetilde {(\R^d)^n}$.
 
 Let us introduce a metric compatible with the described topology on $\Go$.
 We set 
 \begin{equation*}
 dist(\zeta,\eta):=\left\{\begin{array}{ll} 
 1 \wedge d_{Eucl}(\zeta,\eta), &|\zeta|=|\eta|, \\
 1, &\textrm{otherwise}. \end{array}\right.
 \end{equation*}
Here
$d_{Eucl}(\zeta, \eta) $ is the metric induced by the Euclidean metric and 
the map $sym$:
\begin{equation}
 \begin{split}
d_{Eucl}(\zeta, \eta) 
= \inf \{|x-y| : x \in sym ^{-1}  \zeta , y \in sym ^{-1} \eta  \},
\end{split}
\end{equation}
where
$|x-y|$ is the Euclidean distance between $x$ and $y$,
$sym ^{-1} \eta =sym ^{-1} ( \{ \eta \} )$. 
In many aspects, this metric 
resembles the
Wasserstein type distance 
in \cite{RockSchied}. The differences are,
$dist$ is bounded by $1$ and 
it is defined on $\Go$ only.

Note that the metric $dist$ satisfies equalities

\begin{equation}\label{vernachlassigen}
 dist(\zeta \cup x, \eta \cup x ) = dist(\zeta, \eta )
\end{equation}
 for $\zeta , \eta \in  \CYRG _0 (\R^d) $, $x \in \R^d$,
$x \notin \zeta , \eta$,
and 

\begin{equation}\label{zurechtkommen}
 dist(\zeta \setminus x, \eta \setminus x ) = dist(\zeta, \eta ),
\end{equation}
$x \in \zeta, \eta$.
We note that the space $\Go$ equipped with this metric is not 
complete. Nevertheless, $\Go$ is a Polish space, i.e.,
$\Go$ is separable and
there exists a metric $\tilde \rho$ which induces 
the same topology as $dist$ does and such
that $\Go$ equipped with $\tilde \rho$ is a complete metric space. 
To prove this, we embed $\CYRG _0 ^{(n)}(\R ^d)$ into 
the space $\ddot{\CYRG} _0 ^{(n)}(\R ^d)$ of $n$-point
multiple configurations, which 
we define as the space of all counting 
measures $\eta$ on $\R ^d$ with $\eta (\R ^d) = n$.
Abusing notation,
we may represent each 
$\eta \in \ddot{\CYRG} _0 ^{(n)}(\R ^d)$ as a set
$\{x_1,...,x_n \}$, where some points
among $x_j \in \R ^d$
may be equal (recall our convention on
identifying a configuration with a measure; as a 
measure, $\eta = \sum\limits _{j=1} ^n \delta _{x_j}$). 
One should keep in mind that $\{x_1,...,x_n \}$ is
not really a set here, since it is possible that 
$x_i = x_j$ for $i \ne j$, $i,j \in \{1,...,n \}$.
The representation allows us to extend 
$sym$ to the map 

 \begin{equation}
  \begin{split}
  \overline{sym} : \bigsqcup\limits _{m=0}^\infty  {(\R^d)^n} \rightarrow 
  \ddot{\CYRG} _0 ^{(n)}(\R ^d) \\
    \overline{sym}((x_1,...,x_n)) := \{x_1 ,..., x_n  \},
  \end{split}
 \end{equation}
and
define a metric on $\ddot{\CYRG} _0 ^{(n)}(\R ^d)$:
for $\zeta, \eta \in \ddot{\CYRG} _0 ^{(n)}(\R ^d)$
we set $\overline{dist}(\zeta, \eta ) =1 \wedge \overline {d_{Eucl}} (\zeta , \eta) $,
$\overline {d_{Eucl}}(\zeta, \eta) $ is the metric induced by the Euclidean metric and 
the map $\overline{sym}$:

\begin{equation}
 \begin{split}
\overline {d_{Eucl}}(\zeta, \eta) 
= \inf \{|x-y| : x \in \overline{sym} ^{-1}  \zeta , y \in \overline{sym} ^{-1} \eta  \},
\end{split}
\end{equation}
 
 The metrics $dist$ and $\overline{dist}$ coincide on 
 $\CYRG _0 ^{(n)}(\R ^d) \times \CYRG _0 ^{(n)}(\R ^d)$ (as functions).
 Furthermore, one can see that $(\ddot{\CYRG} _0 ^{(n)}(\R ^d),
 \overline{dist})$
 is a complete separable metric space, and thus a Polish space.
 The next lemma describes convergence 
 in $\ddot{\CYRG} _0 ^{(n)}(\R ^d)$ (compare with
 Lemma 3.3 in \cite{KondKut}).
 
 \begin{lem5}\label{caveat}
  Assume that $\eta ^m \to \eta$ in $\ddot{\CYRG} _0 ^{(n)}(\R ^d)$, 
  and let $\eta = \{x_1,...,x_n \}$. Then 
  $\eta ^m$, $m \in \N$, 
  may be numbered, $\eta ^m = \{x_1 ^m,...,x_n ^m\}$, in such a way 
  that
  
  $$
  x_i ^m \to x_i ,  \ \ \ m \to \infty
  $$
  in $\R ^d$.

 \end{lem5}

   \textbf{Proof}. The inequality $\overline{dist} (\eta ^m, \eta ^m) < \varepsilon $
   implies existence of 
   a point from $\eta ^m$ in the ball $B_\varepsilon(x_i)$ 
   of radius $\varepsilon$
    centered at $x_i$, $i \in \{1,...,n \}$.
   Furthermore, in the case when
   $x_i$ is a multiple point, i.e., if $x_j = x_i$
   for some $j \ne i$, then there are at 
   least as many points
   from $\eta ^m$ in $B_\varepsilon(x_i)$ 
   as $\eta (\{x_i \})$. 
   Observe that, for $\varepsilon <
   \frac12 \inf \{ |x-y|: \eta(\{x\}),\eta(\{x\}) \geq 1 \} \wedge 1$,
   we have in the previous sentence ``exactly as many''
   instead of ``at least as many'', because otherwise there
   would not be enough points in $\eta ^m$. 
   The statement of the lemma follows by letting $\varepsilon \to 0$.

 \begin{lem4}
  $\Go$ is a Polish space.
 \end{lem4}

 \textbf{Proof}. Since $\Go$ is a disjoint union of countably 
 many spaces $\CYRG _0 ^{(n)}(\R ^d)$,
 it  suffices to establish that each of
 them is a Polish space. To prove that $\CYRG _0 ^{(n)}(\R ^d)$
 is a Polish space, $n \in \N$, we will show that
 it is a countable intersection of open sets
 in a Polish space $\ddot{\CYRG} _0 ^{(n)}(\R ^d)$.
 Then we may apply Alexandrov's theorem:
 any $G_{\delta}$ subset of a Polish space is a Polish
 space, see  \S 33, VI in \cite{Kuratowski}.
 
 To do so, denote by $\mathbf{B}_m$ the closed
 ball of radius $m$
  in $\R ^d$, with the center at the origin.
 Define $F_m : =\{ \eta \in \ddot {\CYRG} _0 ^{(n)}(\R ^d) \mid
 \eta (\{x \}) \geq 2 \textit{ for some } x \in \mathbf{B} _m \}$
 and note that 
 
 $$
  {\CYRG} _0 ^{(n)}(\R ^d) = \bigcap _{m=1} ^ \infty  [\ddot {\CYRG} _0 ^{(n)}(\R ^d) 
  \setminus F_m]
 $$
 
 Since $\ddot {\CYRG} _0 ^{(n)}(\R ^d)$ is Polish, it only remains 
 to show that $F_m$ is closed in $\ddot {\CYRG} _0 ^{(n)}(\R ^d)$.
 This is an immediate consequence of the previous lemma.

 \subsubsection{Lebesgue-Poisson measures}\label{Leb-Pos measure}
 
 Here we define the Lebesgue-Poisson measure on $\Go$,
 corresponding to a non-atomic Radon measure $\sigma$ on $\R ^d$. 
 Our prime example for $\sigma$ will be the Lebesgue measure on $\R ^d$. 
 For any $n \in \N$ the product measure
 $\sigma ^{\otimes n} $ can be considered
 by restriction as a measure on $\widetilde {(\R^d)^n}$.
 The projection of this measure on 
 $\CYRG _0 ^{(n)}$ via $sym$ we denote by 
 $\sigma ^{ (n)} $, so that
 
 \[
  \sigma ^{ (n)} (A) = \sigma ^{\otimes n} (sym^{-1} A), \ \ \ \ \ A \in \mathscr{B}(\CYRG _0 ^{(n)}).
 \]

 On $\CYRG _0 ^{(0)}$
 the measure $\sigma ^{ (0)} $ is given by 
 $\sigma ^{ (0)}(\{ \varnothing \}) =1$.
 The \textit{Lebesgue-Poisson measure} on $(\Go, \mathscr{B}(\Go))$
 is defined as 
 
 \begin{equation}
  \lambda _{\sigma} := \sum\limits _{n=0} ^\infty \frac{1}{n!} \sigma ^{ (n)}.
 \end{equation}

 The
measure $\lambda _{\sigma}$ is finite iff $\sigma$ is finite. 
We say that $\sigma$ is the \emph{intensity measure} of $\lambda _{\sigma}$.

 \subsubsection{The Skorokhod space}\label{D space}
 
 For a complete separable metric space $(E,\rho)$ the space $D_E$ of 
 all cadlag $E$-valued functions equipped 
 with the Skorokhod topology is a Polish space; 
 for this statement and related definitions, see, e.g.,
 Theorem 5.6, Chapter 3 in \cite{EthierKurtz}.
 Let $\rho _D$ be a metric on $D_E$
 compatible with the Skorokhod topology and such that
 $(D_E , \rho _D)$ is a complete separable metric space.
  Denote by $(\mathcal P (D_E),\rho _p)$ the metric space
  of probability measures on $\mathscr{B}(D_E)$,
  the Borel $\sigma$ - algebra of $D_E$,
  with the Prohorov metric, i.e. for $P,Q \in \mathcal P (D_E)$
  
  \begin{equation}
   \rho _p (P,Q) = \inf \{ \varepsilon >0 : P(F)\leq 
   Q(F^{\varepsilon}) + \varepsilon  \text{ for all } F \in \mathscr{B}(D_E) \}   
  \end{equation}
where
$$F^{\varepsilon} = \{x\in D_E: \rho _D (x,F) < \varepsilon  \}.$$
 
 Then $(\mathcal P (D_E),\rho _p)$ is separable and complete;
 see, e.g., \cite{EthierKurtz}, Section 1, Chapter 3, and
 Theorem 1.7, Chapter 3.
 The Borel $\sigma$-algebra $\mathscr{B}(D_E)$ 
 coincides with the one generated by the coordinate mappings;
 see Theorem 7.1, Chapter 3 in \cite{EthierKurtz}.
 In this work, we mostly consider $D_{\Go} [0;T]$ and $D_{\G} [0;T]$ endowed with
 the Skorokhod topology.

 \subsection{Integration with respect to Poisson point processes}

We give a short introduction to the theory of
integration with respect to Poisson point processes.
For construction of Poisson point processes
with given intensity, see e.g.
\cite[Chapter 12]{KallenbergFound},
\cite{KingmanPP}, \cite[Chapter 12, $\mathsection$ 1]{RevuzYor}
 or \cite[Chapter 1, $\mathsection$ 8,9]{IkedaWat}.
All definitions,
constructions and statements about 
integration given here may 
be found in \cite[Chapter 2, $\mathsection$ 3]{IkedaWat}.
See also \cite[Chapter 1]{GikhSkor3}
for the theory of integration 
with respect to an orthogonal
martingale measure.

On some filtered
 probability space $(\Omega , \mathscr{F}, \{ \mathscr{F} \}_{t\geq 0}, P)$,
consider a Poisson point process $N$ on
$\R _+ \times \X \times \R _+$ with
intensity measure $dt\times \beta (dx) \times du$,
where $\X = \R ^d$
or $\X = \Z ^d$. 
We require the filtration $\{ \mathscr{F} \}_{t\geq 0}$
to be increasing and right-continuous, and 
we assume that $\mathscr{F}_0$ is complete 
under $P$.
We interpret the argument from the first space
$\R _+$ as time. For
$\X=\R ^d$ the intensity
measure $\beta$ will be the Lebesgue
measure on $\R ^d$, for
 $\X = \Z ^d$ we set $\beta = \#$, where
 \[
  \# A = |A|, \ \ \  A \in \mathscr{B}(\Z ^d).
 \]
The Borel $\sigma$-algebra over
$\Z ^d$ is the collection 
of all subsets of $\Z ^d$, i.e.
$\mathscr{B}(\Z ^d) = 2^{\Z ^d}$.
Again, as is the case with configurations, for $X = \R ^d$
we treat a point process as a random collection 
of points as well as a random measure.

We say that the process
$N$ is called \textit{compatible} with 
 $(\mathscr{F}_t, t \geq 0)$
if $N$ is adapted, that is, all random variables
of the type $N(\bar T _1, U)$, $\bar T _1 \in \mathscr{B}([0;t])$,
$U \in \mathscr{B}(\X \times \R _+)$, are $\mathscr{F}_t$-measurable,
and all random variables of the type
$N(t+h,U) - N(t,U)$, $ h \geq 0 , U \in \mathscr{B}(\X \times \R _+)$,
are independent of $\mathscr{F}_t$, $N(t,U) = N([0;t], U)$.
For any $U \in \mathscr{B}(\X \times \R _+)$ with
$(\beta \times l) (U) < \infty$, $l$ is the
Lebesgue measure on $\R ^d$,
the process $(N([0;t], U) - t \beta \times l (U), t\geq 0)$
is a martingale (with respect to $(\mathscr{F}_t, t \geq 0)$;
see \cite[Lemma 3.1, Page 60]{IkedaWat}). 

\begin{def5}

A process $f: \R _+ \times \X \times \R _+ \times \Omega \to \R $ is 
\textit{predictable}, if it is measurable with
respect to the smallest $\sigma$ - algebra
generated by all $g$ having the following properties:

(i) for each $t>0$, $(x,u, \omega) \mapsto g(t, x, u , \omega))$
is $\mathscr{B}(\X \times \R_+) \times \mathscr{F} _t$ measurable;

(ii) for each $(x,u, \omega)$, the map 
$t \mapsto g(t, x, u , \omega))$ is left continuous.
\end{def5}

For a predictable process 
$f \in L^1 ([0;T] \times \X \times \R _+ \times \Omega)$, 
$t \in [0;T]$ and $U \in \mathscr{B}(\X \times \R _+)$ we define 
the integral 
$
 I_t(f)= \int\limits _{[0;t] \times U} f(s, x, u , \omega)dN(s,x,u)
$
as the Lebesgue-Stieltjes integral
with respect to the measure $N$:

\[
 \int\limits _{[0;t] \times U} f(s, x, u , \omega)dN(s,x,u) = 
 \sum\limits _{s \leq t, (s,x,u) \in N} f(s, x, u , \omega).
\]

This sum is well defined, since 

\[
 E \sum\limits _{s \leq t, (s,x,u) \in N} |f(s, x, u , \omega)|=
 \int\limits _{[0;t] \times U} |f(s, x, u , \omega)|ds \beta (dx) du < \infty
\]
We use $dN(s,x,u)$ and $N(ds,dx,du)$ interchangeably when
we integrate over all variables.
The process $I_t(f)$ is right-continuous as a function of $t$,
and adapted. Moreover, the process

\[
 \tilde I_t(f)  = \int\limits _{[0;t] \times U} f(s, x, u , \omega)[dN(s,x,u) - ds \beta (dx) du]
\]

is a martingale with respect to $(\mathscr{F}_t, t \geq 0)$, \cite[Page 62]{IkedaWat}. Thus,

\begin{equation}
 E \int\limits _{[0;t] \times U} f(s, x, u , \omega)dN(s,x,u) = 
 E \int\limits _{[0;t] \times U} f(s, x, u , \omega) ds \beta (dx) du.
\end{equation}

This equality will be used several times throughout this work.

\begin{rmk9}

We can extend the collection of integrands, in particular,
we can define \\
$\int\limits _{[0;t] \times U} f(s, x, u , \omega)dN(s,x,u)$
for $f$ satisfying

\[
 E\int\limits _{[0;t] \times U} (|f(s, x, u , \omega)| \wedge 1) ds \beta (dx) du < \infty.
\]

However, we do not use such integrands.
\end{rmk9}

The Lebesgue-Stieltjes integral is
defined $\omega$-wisely and it is
a function
of an integrand and an integrator. As a result, we have
 the following statement. The
sign $\,{\buildrel d \over =}\,$ means
equality in distribution.

\begin{3stm} \label{distributions}
Let $M _k$ be Poisson point processes defined on some, possibly different,
probability spaces, and let
$\alpha _k$ be integrands, $k=1,2$,
such that integrals $\int \alpha_k  dM _k$
are well defined. If 
$(\alpha _1 , M _1 ) \,{\buildrel d \over =}\, (\alpha _2 , M _2 )$,
then 
\[
\int \alpha_1  dM _1 \,{\buildrel d \over =}\, \int \alpha_2  dM _2.
\]
\end{3stm}

The proof is straightforward.

\subsubsection{An auxiliary construction}\label{voranbringen}

Let $\tilde{\#}$ be the counting measure on $[0,1]$, i.e. 
\[
 \tilde \# C  = |C|, \ \ \ \ \ C \in \mathscr{B}([0;1]).
\]
The 
measure $\tilde \#$ is not $\sigma$-finite. For a 
cadlag $\Go$-valued process $(\eta _t)_{t \in [0;\infty]}$,
adapted to $\{\mathscr{F}_t \}_{t \in [0;\infty]}$,
we would like to define integrals of the form

\begin{equation}\label{sublime}
 \int\limits _{\R ^d \times [0; \infty ]\times [0;\infty) } I _{\{x\in B \cap \eta _{r-}\}}
f(x, r,v, \omega )   d \tilde N_2 (x,r,v) 
\end{equation}
where $B$ is a bounded Borel subset of $\R^d$,
$f$ is a bounded predictable process and $ \tilde N_2$
is a Poisson point process  on 
$\R ^d \times [0;T]\times [0;\infty)$ with intensity $ \tilde \# \times dr \times dv$,
compatible with $\{\mathscr{F}_t \}_{t \in [0;\infty]}$.

We can not hope to give a meaningful definition 
for an integral of the type \eqref{sublime},
because of the measurability issues. For example, 
the map
\begin{align*}
 \Omega &\to \R, \\
 \omega &\mapsto \tilde N_2  (u(\omega),[0;1],[0;1]),
\end{align*}
where 
$u$ is an independent of $ \tilde N_2$
uniformly distributed on $[0;1]$
random variable, does not have
to be a random variable.
Even if it were a random variable, some undesirable 
phenomena would appear, see, e.g., \cite{Uncountable}.

To avoid this difficulty, we employ another construction. A
similar approach was used in \cite{FournierMeleard}.
If we could give meaningful definition to the
integrals of the type \eqref{sublime}, we would 
expect 
\begin{align*}
 \int\limits _{\R ^d \times [0; t ]\times [0;\infty) } I& _{\{x\in B \cap \eta _{r-}\}}
f(x, r,v, \omega )   d \tilde N_2 (x,r,v) - \\
 \int\limits _{\R ^d \times [0; t ]\times [0;\infty) } &I _{\{x\in B \cap \eta _{r-}\}}
f(x, r,v, \omega )  \tilde \#(dx)drdv
\end{align*}
to be a martingale (under some conditions on $f$ and $B$).

Having this in mind, consider a  Poisson point process   $  N_2$ on 
$ \Z  \times \R _+ \times \R _+$ with intensity $\# \times dr \times dv$,
defined on $(\Omega , \mathscr{F}, \{ \mathscr{F} \}_{t\geq 0}, P)$
(here $\#$ denotes the counting measure on $\Z$. This measure 
is $\sigma$-finite). 
We require  $ N _2$ to be compatible with
 $\{ \mathscr{F} \}_{t\geq 0}$.
Let $(\eta_t)_{t\in[0,\infty]}$ be an adapted cadlag process in $\Gamma_0(\mathbb{R}^d)$, 
satisfying the following condition: for any $T<\infty$, 

\begin{equation}\label{back to back}
 R_T = |\bigcup _{t \in [0;T]} \eta _t| < \infty \ \ \ \  \textrm{a.s.}
\end{equation}
The 
set $R_{\infty}: = \bigcup _{t \in [0;\infty]} \eta _t$ is at most countable, provided
\eqref{back to back}.
Let $\preccurlyeq$ be the lexicographical order on $\R^d$. We 
can label the points of $\eta _0$, 

\[
\eta _0 = \{ x_{0},x_{-1},...,x_{-q}  \}, \ \ \ x_{0} \preccurlyeq x_{-1} \preccurlyeq  ... \preccurlyeq x_{-q}.
\]
There exists an a.s. unique representation
\[
R _{\infty}  \setminus \eta _0 = \{x_1,x_2,... \}
\]
such that for any $n,m \in \N$, $n <m$,
either $\inf\limits _{s\geq 0} \{s: x_n \in \eta _s  \} < \inf\limits _{s\geq 0} \{s: x_m \in \eta _s  \}$,
or $\inf\limits _{s\geq 0} \{s: x_n \in \eta _s  \} = \inf\limits _{s\geq 0} \{s: x_m \in \eta _s  \}$
and $x_n \preccurlyeq x_m$. In other words, as time goes on,
appearing points are added to $\{x_1,x_2,... \}$ in the order in 
which they appear. If several points appear 
simultaneously, we add them in the lexicographical order.

For the sake of convenience, we set $x_{-i} = \Delta$, $i \leq -q-1$,
where $\Delta \notin \Z$.
We say that the sequence 
$\{...,x_{-1},x_1,x_2,... \}$ is \emph{related} to $(\eta _t)_{t \in [0;\infty]}$.

For a predictable process 
$f \in L^1 (\R ^d \times \R _+ \times \R _+ \times \Omega)$ and $B\in \mathscr{B}(\R ^d)$, 
consider

\begin{equation}
 \int\limits _{\Z \times (t_1;t_2] \times [0; \infty ) } I_{\{x_i \in \eta _{r-} \cap B \}}
f(x_i, r,v, \omega )   d  N _2 (i,r,v). 
\end{equation}

Assume that $R_T$ is bounded for some $T>0$.
Then, for a bounded predictable $f \in L^1 (R^d \times \R _+ \times \R _+ \times \Omega)$
and $B\in \mathscr{B} (\R ^d)$,
the process 
\[
\int\limits _{\Z \times (0;t] \times [0; \infty ) }  I_{\{x_i \in \eta _{r-} \cap B\}}
f(x_i, r,v, \omega )  d  N_2 (i,r,v)
\]
\[
- \int\limits _{\Z \times (0;t] \times [0; \infty ) } I_{\{x_i \in \eta _{r-} \cap B\}}
f(x_i, r,v, \omega )  \#(di)drdv
\]
is a martingale, cf. \cite[Page 62]{IkedaWat}.

\subsubsection{The strong Markov property of a Poisson point process}

We will need the strong Markov property of a Poisson point process.
To simplify notations, assume that $N$ is a Poisson point process on
$\R _+ \times \R ^d $ with
intensity measure $dt\times dx$. 
Let $N$ be compatible with a right-continuous complete filtration
$\{ \mathscr{F}_t \}_{t \geq 0}$, and $\tau$ be a finite a.s.
$\{ \mathscr{F}_t \}_{t \geq 0}$-stopping time 
 (stopping time with respect to $\{ \mathscr{F}_t \}_{t \geq 0}$). Introduce another
Point process $\overline N $ on $\R _+ \times \R ^d $,

\[
 \overline N ([0;s] \times U) = N ((\tau;\tau + s] \times U), \ \ \  U \in \mathscr{B}(\R ^d).
\]

\begin{prop1}\label{Strong MP}
 The process $\overline N$ is a Poisson point process with
 intensity $dt\times dx$, independent of 
 $\mathscr{F}_{\tau}$.
\end{prop1}

\textbf{Proof}. To prove the proposition, it is enough to show that

(i) for any $b>a>0$ and open bounded $U \subset \R ^d$, 
$\overline N ((a;b),U)$ is a Poisson random variable
with mean $(b-a) \beta (U)$, and

(ii) for any $b_k>a_k>0$, $k=1,...,m$, and any open bounded
$U_k \subset \R ^d$, such that $((a_i;b_i) \times U_i) \cap( (a_j;b_j) \times U_j) = \varnothing $, $i \ne j$,
 the collection $\{\overline N ((a_k;b_k) \times U_k)\}_{k=1,m}$ 
 is a sequence of independent random variables, independent of $\mathscr{F}_{\tau}$.

Indeed, $\overline N $ is determined completely by values on sets of type
$(b-a) \beta (U)$, $a,b,U$ as in (i), therefore it must be
an independent of $\mathscr{F}_{\tau}$ Poisson
point process 
if (i) and (ii) hold.

Let $\tau _n$ be the sequence of $\{ \mathscr{F}_t \}_{t \geq 0}$-stopping times, 
$\tau _n = \frac{k}{2^n}$ on $\{ \tau \in (\frac{k-1}{2^n};\frac{k}{2^n}] \}$, 
$k \in \N$. Then
 $\tau _n \downarrow \tau$ and
$\tau _n - \tau \leq \frac{1}{2^n}$. The
stopping times $\tau _n$ take only countably many values.
The process $N$ satisfies the strong Markov property for $\tau _n$: 
 the processes $\overline N _n$, defined by
$$
\overline N _n ([0;s] \times U) := N ((\tau_n;\tau_n + s] \times U), 
$$
are Poisson point processes, independent of $\mathscr{F}_{\tau _n}$.
To prove this, take $k$ with $P\{ \tau_n = \frac{k}{2^n} \} >0$ and
note that on $\{ \tau_n = \frac{k}{2^n} \}$,  $\overline N _n $ 
coincides with
 process the Poisson point process $\tilde N _{\frac{k}{2^n}}$
given by
\[
\tilde N _{\frac{k}{2^n}} ([0;  s] \times U)
:=  N  \bigg((\frac{k}{2^n};\frac{k}{2^n} + s] \times U) \bigg), \ \ \  U \in \mathscr{B}(\R ^d).
\]
Conditionally on $\{ \tau_n = \frac{k}{2^n} \}$,
$\tilde N _{\frac{k}{2^n}}$ is again a Poisson point process,
with the same intensity. Furthermore, conditionally on $\{ \tau_n = \frac{k}{2^n} \}$,
$\tilde N _{\frac{k}{2^n}}$
is independent of $\mathscr{F}_{\frac{k}{2^n}}$, hence it
is independent 
of 
$ \mathscr{F}_{\tau } \subset \mathscr{F}_{\frac{k}{2^n}}$.

To prove (i), 
note that $\overline N _n ((a;b) \times U) \to \overline N ((a;b) \times U) $ a.s. 
and all random variables $\overline N _n ((a;b) \times U)$ have the same distribution, 
therefore $\overline N  ((a;b) \times U)$ is a Poisson random variable with mean 
$(b-a)\lambda (U)$. The random variables $\overline N _n ((a;b) \times U)$
are independent of $\mathscr{F}_{\tau}$, 
hence $\overline N  ((a;b) \times U)$ is independent 
of $\mathscr{F}_{\tau}$, too. Similarly, (ii) follows. $\Box$

Analogously, the strong Markov property for 
a Poisson point process
on $ \R_+ \times \N$
with intensity $dt \times \#$
may be formulated and proven.

\begin{rmk7}\label{Missstand}
We assumed in Proposition \ref{Strong MP}
that the filtration $\{ \mathscr{F}_t \}_{t \geq 0}$,
compatible with $N$,
is right-continuous and complete. To be able to apply 
Proposition \ref{Strong MP},
we should show that such filtrations exist.

Introduce the natural filtration of $N$,

\[
 \mathscr{F}_t ^0 = \sigma \{  N_k(C, B), 
B\in \mathscr{B} (\R ^d), C\in \mathscr{B} ([0;t]) \},
\]
and let $\mathscr{F}_t $ be the completion of $\mathscr{F}_t ^0$
under $P$.
Then $N$ is compatible with $\{\mathscr{F}_t\} $.
We claim that $\{ \mathscr{F}_t \}_{t \geq 0}$, defined in such a way,
is right-continuous (this may be regarded as an 
analog of Blumenthal $0-1$ law). Indeed, as in the proof 
of Proposition \ref{Strong MP}, one may check
that $\tilde N _a$ is independent of $\mathscr{F}_{a+}$. Since 
$\mathscr{F}_{\infty}  = \sigma(\tilde N _a) \vee \mathscr{F}_{a}$,
$\sigma(\tilde N _a)$ and $\mathscr{F}_{a}$ are independent
and $\mathscr{F}_{a+} \subset \mathscr{F}_{\infty} $, one sees that
$\mathscr{F}_{a+} \subset \mathscr{F}_{a} $. Thus, 
$\mathscr{F}_{a+} = \mathscr{F}_{a} $.

\end{rmk7}

\begin{rmk8}
 We prefer to work with right-continuous complete filtrations,
 because we want to
  ensure that there is no problem
 with conditional probabilities, and that the
 hitting times we will consider are stopping times.
\end{rmk8}

  \subsection{Miscellaneous}
  
   When we write $\xi  \sim Exp( \lambda)$, we mean
  that the random variable $\xi$ is
 exponentially distributed 
 with parameter $\lambda$.
  
  \begin{lem7}\label{preclude}
   If $\alpha$ and $\beta$ are exponentially distributed random variables
   with parameters $a$ and $b$ respectively 
   (notation: $\alpha \sim Exp(a)$, $\beta \sim Exp(b)$ )
   and they are independent, then
   
   \[
   P\{ \alpha < \beta   \} = \frac{a}{a+b}.
   \]
   
  \end{lem7}
  Indeed, 
  
  \[
   P\{ \alpha < \beta   \} = \int _0 ^\infty a P\{ x < \beta   \}e^{-ax}
   = a \int _0 ^\infty e^{-(a+b)x} = \frac{a}{a+b}.
   \]

  Here are few other properties of exponential 
  distributions.
  If $\xi _1,\xi _2, ..., \xi _n$
  are independent exponentially distributed 
  random variables with parameters
  $c_1,...,c_n$ respectively, then 
  $\min\limits_{k\in \{1,...,n \}} \xi _k$
  is exponentially distributed
  with parameter $c_1+...+c_n$.
  Again, the proof may be done
  by direct computation.
  If $\xi _1,\xi _2, ...$
  are independent exponentially distributed 
  random variables with parameter
  $c$ and $\alpha _1,\alpha _2, ...$
  is an independent sequence of independent
  Bernoulli 
  random variables with parameter $p \in (0;1)$,
  then the random variable
  
  \[
   \xi = \sum\limits _{i=1}^\theta \xi _i, \ \ \
   \theta= \min \{ k\in \N : \alpha _k =1 \}
  \]
is exponentially distributed with parameter $\frac{c}{p}$.
  The random variable $\xi$
  is the time of the first jump 
  of a thinned Poisson point process with
  intensity $c$.
  The statement about the 
  distribution of $\xi$ is a consequence of the property that the independent
  thinning of a Poisson point process with
  intensity $\lambda$ is a 
  Poisson point process 
  with intensity $p\lambda$, see \cite[Theorem 12.2,(iv)]{KallenbergFound}.

  We will also need the result about finiteness of
  the expectation of the Yule process. A Yule
  process $(Z_t)_{t\geq 0}$ is a  pure birth Markov process in $\Z _+$
  with birth rate $\mu n$, $\mu >0$, $n \in \Z_+$. That is,
  if $Z_t = n$, then a birth occur 
  at rate $\mu n$, i.e.
  \[
   P\{Z_{t+\Delta t} - Z_t = 1 \mid Z_t = n \} = \mu n + o(\Delta t).
  \]
For more details \label{Yule page}
  about Yule processes see e.g. \cite[Chapter 3]{Branch1},
  \cite[Chapter 5]{Branch2}, \cite{Yuleproc}
  and references therein. Let $(Z_t (n))_{t\geq 0}$
  be a Yule process started at $n$. 
  The process $(Z_t (n))_{t\geq 0}$ can 
  be considered 
  as a sum of $n$ independent Yule processes started from $1$, see e.g. \cite{Yuleproc}.
  The expectation of $Z_t(1) $ is finite 
  and $ E Z_t(1) = e^{\mu t}$, see e.g.
  \cite[Chapter 3, Section 6]{Branch1} or 
  \cite[Chapter 5, Sections 6,7]{Branch2}.
  Consequently, if $(Z_t)_{t\geq 0}$
  is a Yule process with $E Z_0 < \infty$,
  then $E Z_t < \infty$ and $E Z_t = E Z_0 e^{\mu t}$.
  
  Here are some other properties of 
  Poisson point processes which are used throughout in the article.
  If $N$ is a Poisson point process on $\R _+ \times \R ^d \times \R _+$
  with intensity $ds \times dx \times du$, then a.s.
  \begin{equation} \label{travesty}
   \forall x \in \R ^d : N(\R _+ \times \{x \} \times \R _+) \leq 1.
\end{equation}
  Put differently, no plane of the form $\R _+ \times \{x \} \times \R _+$
  contains more than $1$ point of $N$. Using the $\sigma$-additivity of 
  the probability measure,
  one can deduce \eqref{travesty}
  from
   \begin{equation} \label{sentient}
   \forall x \in \R ^d : N([0;1] \times \{x \} \times [0;1]) \leq 1.
\end{equation}
  We can write 
  \[
   \bigg\{\forall x \in \R ^d : N([0;1] \times \{x \} \times [0;1]) \leq 1 \bigg\} 
   \]
   \[
   \supset
   \bigg\{ \forall k \in \{0,1,...,n-1 \} :
   N([0;1] \times [\frac{k}{n};\frac{k+1}{n}] \times [0;1]) \leq 1 \bigg\},
  \]
 and then we can compute
 \[
 P \bigg\{ \forall k \in \{0,1,...,n-1 \} :
   N([0;1] \times [\frac{k}{n};\frac{k+1}{n}] \times [0;1]) \leq 1 \bigg\} 
   \]
   \[= 
  \Big( P \{
   N([0;1] \times [0;\frac{1}{n}] \times [0;1]) \leq 1 \} \Big)^{n}=
  \Big( \exp (-\frac 1n)[1 + \frac 1n] \Big)^{n} = \Big(1 - o(\frac 1n)\Big)^{n} =1 -o(\frac 1n).
 \]
Thus, \eqref{sentient} holds.
  
  Let $\psi \in L ^1 (\R ^d)$, $\psi \geq 0$. Consider the time until
  the first arrival 
  
\begin{equation}
   \tau = \inf \{t>0: \int\limits _{[0;t] \times \R ^d \times \R _+ } I_{[0;\psi(x)]} (u) N(ds,dx,du) >0\}.
\end{equation}
The random variable $\tau$ is distributed exponentially with the parameter
$||\psi||_{L^1}$. From \eqref{travesty} we know that a.s.
\[
 N(\{ \tau \} \times \R ^d  \times \R _+) =  
 N \bigl( \{(\tau,x,u) \mid x\in \R ^d, u \in [0;\psi(x)] \} \bigr) =1
\]
Let $x_{\tau}$ be the unique element of $\R ^d$ defined by
\[
 N(\{ \tau \} \times \{ x_{\tau} \}   \times \R _+) =1.
\]
Then
\begin{equation}
 P\{x_{\tau} \in B \} = \frac{\int _B \psi (x) dx }{\int _{\R ^d} \psi (x) dx}, \ \ \ 
 B \in \mathscr{B} (\R ^d).
\end{equation}

\subsection{Pure jump type Markov processes}\label{Pure jump type MP}

In this section we give a very concise treatment of 
pure jump type Markov processes. Most of 
the definitions and facts given here
can be found in 
\cite[Chapter 12]{KallenbergFound};
see also, e.g., \cite[Chapter 3, $\mathsection$ 1]{GikhSkor2}.

We say that a process $X = (X_t)_{t \geq 0}$ in some measurable 
space $(S , \mathcal{S})$ 
is of \textit{pure jump type} if 
 its paths are a.s. right-continuous and constant apart from
isolated jumps. In that case we may denote the jump times 
of $X$ by $\tau _1,\tau _2,...$, with understanding that
$\tau _n = \infty$ if there are fewer that $n$ jumps.
The times $\tau _n$ are stopping times with 
respect to the right-continuous filtration induced by $X$.
For convenience we may choose $X$ to be the identity mapping 
on the canonical path space $(\Omega,\mathscr{F})= 
(S^{[0;\infty)} , \mathcal{S} ^{[0;\infty)})$.
When $X$ is a Markov process,
the distribution with initial state $x$ is denoted 
by $P_x$, and we note that the mapping $x \mapsto P_x(A)$ 
is measurable in $x$, $A \in \Omega$.

\textbf{Theorem 12.14} \cite{KallenbergFound} 
\textit{(strong Markov property, Doob) A pure jump 
type Markov process satisfies strong Markov property at every
stopping time.}

We say that a state $x\in S$ is \emph{absorbing}
if $P_x \{ X \equiv x \} =1$.

\textbf{Lemma 12.16} \cite{KallenbergFound}
\textit{If $x$ is non-absorbing,
then under $P_x$ the time $\tau_1$ until the first jump
is exponentially distributed and independent
of $\theta _{\tau _1}X$.
}

Here $\theta _t$ is a shift, and $\theta _{\tau _1}X$ defines a 
new process, 

\[
 \theta _{\tau _1}X (s)=  X(s+ \tau _1).
\]

For a non-absorbing state $x$, we may define the \textit{rate function}
$c(x)$ and \textit{jump transition kernel} 
$\mu (x,B)$ by

\[
 c(x) = (E_x \tau _1 )^{-1} , \ \
 \mu (x,B) = P _x \{ X_{\tau _1} \in B \} , 
 \ \  \ \    x\in S , \ B \in \mathcal{S}.
\]

In the sequel, $c(x)$ will also be referred to as \textit{jump rate}.
The kernel $c \mu $ is called a 
\textit{rate kernel}.

The following theorem gives an explicit representation
of the process in terms of a discrete-time Markov chain 
and a sequence of exponentially distributed
random variables. This result shows in particular
that the distribution $P_x$ is uniquely determined by the 
rate kernel $c \mu $. We assume existence of the required
randomization variables (so that the underlying
probability space is ``rich enough'').

\textbf{Theorem 12.17} \cite{KallenbergFound} \label{thm 12.17 Kall}
\textit{(embedded Markov chain) Let $X$ be 
a pure jump type Markov process with rate 
kernel $c \mu $. Then there exists a Markov process
$Y$ on $\Z _+$ with transition 
kernel $\mu$ and an independent sequence of i.i.d.,
exponentially distributed random variables
$\gamma _1, \gamma _2, ...$ with mean $1$ such that
a.s.}

\begin{equation}\label{beforehand}
 X_t = Y _ n , \ \ \ t \in [\tau _n, \tau _{n+1}), \  n \in \Z_+ ,
\end{equation}
\textit{where}
\begin{equation}\label{ruby}
 \tau _n = \sum\limits _{k=1} ^n \frac{\gamma _k}{c(Y_{k-1})}, 
 \ \ n \in \Z_+.
\end{equation}

In particular, the differences between the moments of jumps 
$ \tau _{n+1} - \tau _n$
of a pure \label{Einoede}
jump type Markov process are exponentially distributed
given the embedded chain $Y$, with parameter $c(Y_n)$.
If $c(Y_{k}) =0$
for some (random) $k$, we set
$\tau _n = \infty$ for $n \geq k+1$, while
 $Y_{n}$ are not defined, $n \geq k+1$.

\textbf{Theorem 12.18} \cite{KallenbergFound} \label{synthesis}
\textit{(synthesis) For any rate kernel $c \mu$ on 
$S$ with $\mu(x,\{x \}) \equiv 0$,
consider a Markov chain $Y$ 
with transition kernel $ \mu$ and a
sequence $\gamma _1 ,\gamma _2,... $ of independent
exponentially distributed random
variables with mean 1, independent of $Y$. Assume that 
$\sum _n \frac{\gamma _n}{c(Y_n)} = \infty$
a.s. for every initial distribution
for $Y$. Then \eqref{beforehand}
and \eqref{ruby}
define a pure jump type Markov process
with rate kernel $c \mu$.
}

Next proposition gives a convenient
criterion for non-explosion.

\textbf{Proposition 12.19} \cite{KallenbergFound} 
\textit{(explosion) For any rate kernel $c \mu$
and initial state $x$, let $(Y_n)$
and $(\tau _n) $ be such as in Theorem 12.17. Then a.s.}

\begin{equation}
\tau _n \to \infty \ \ \ \
\textit{iff} \ \ \ \ 
\sum _n \frac{1}{c(Y_n)} = \infty.
\end{equation}
\textit{In particular, $\tau _n \to \infty$ a.s. when 
$x$ is recurrent for $(Y_n)$.}

  \subsection{Markovian functions of a Markov chain}
  
  Let $(S, \mathscr{B}(S))$ be a Polish (state) space. Consider a
  (homogeneous) Markov 
  chain on $(S, \mathscr{B}(S))$ as a family of probability measures
  on $S^ \infty$. Namely,  on the measurable space
  $
  {(\Omega,\mathscr{F}) = (S^ \infty , \mathscr{B}(S ^\infty ))}
 $
  consider a family of probability measures $\{P_s \}_{s \in S}$ 
  such that for the coordinate mappings 
  
  \begin{align*}
 X_n: \Omega &\rightarrow S, \\
 X_n (s_1,s_2,&...)  = s_n
\end{align*}
  
  the process $X = \{X_n \}_{n \in \Z _+ }$ is a Markov chain, 
  and for all $s \in S$

  $$
  P_s \{X_0 =s \} =1,
  $$
  $$
  P_s \{ X_{n+m_j}\in A_j, j=1,...,k_1 \mid \mathscr{F} _n \} 
  = P_{X_n} \{ X_{m_j} \in A_j, j=1,...,k_1  \}.
  $$
  
  Here $A_j \in \mathscr{B} (S)$,
  $m_j \in \N$, $ k_1 \in \N$, $\mathscr{F} _n = \sigma \{ X_1,...,X_n \}$.
 The space $S$ is separable, hence
 there exists a transition probability 
 kernel
 $Q: S \times \mathscr{B} (S) \rightarrow [0;1]$ such that
 
 $$
 Q(s,A) = P_s \{ X_1 \in A \}, \ \ \ s\in S, \ A \in \mathscr{B} (S).
 $$
  
  Consider a transformation of the chain $X$, $Y_n = f(X_n)$, where
  $f:S\to \Z _+$ is a Borel-measurable function,
  with convention $\mathscr{B}(\Z _+) = 2^{\Z _+}$. In the future
  we will need to know 
  when the process $Y = \{Y_n \}_{\Z _+ }$ is 
  a Markov chain. A similar question 
  appeared for the first time in \cite{BRosenblatt}.

 A sufficient condition for $Y $ to be
  a Markov chain is given in 
  the next lemma.

 \begin{lem6} \label{lumpability}
  Assume that for any bounded Borel function $h: S\rightarrow S$
  
  \begin{equation}\label{insurgent}
   E_s h (X_1) =  E_q h (X_1) \text{\ whenever \ } f(s)= f(q),
  \end{equation}
  Then $Y$ is a Markov chain.
 \end{lem6}
 
 \textbf{Remark}. Condition \eqref{insurgent} is the
 equality of distributions of $X_1$
 under two different measures, $P_s$ and $ P_q$.

 \textbf{Proof}. For the natural filtrations of the
 processes $X$ and $Y$ we have an inclusion
 
 \begin{equation} \label{penultimate}
 \mathscr{F} ^X _n \supset \mathscr{F} ^Y _n, \ \ \  n \in \N,
 \end{equation}
since 
$Y$ is a function of $X$. For $k\in \N$
 and bounded Borel functions $h_j: \Z _+ \rightarrow \R$,
 $j=1,2,...,k$
 (any function on $\Z _+$ is a Borel function), 
 
 \begin{equation} \label{prance}
 \begin{split}
     E_s \left[ \prod\limits _{j=1} ^k 
 h_j(Y_{n+j}) \mid \mathscr{F} ^ X _n  \right] = 
 E_{X_n}  &\prod\limits _{j=1} ^k 
 h_j(f(X_{j})) = \\
  \int _S Q(x_0 , dx_1) h_1(f(x_1)) \int _S Q(x_1 , dx_2) h_2(f(x_2))...
 &\int _S Q(x_{n-1} , dx_n) h_n(f(x_n)) \Bigg| _{x_0 = X _n}
 \end{split}
 \end{equation}

To transform the last integral, we introduce a new kernel:
for $y \in f(S)$ chose $x \in S$ with $f(x) = y$,
ans then for $B \subset \Z _+$
define
\begin{equation}\label{secretive}
 \overline Q (y, B) = Q (x, f^{-1}(B)); 
\end{equation}
The expression on the right-hand side 
 does not depend on the choice of $x$ because of
 \eqref{insurgent}.
 To make the kernel $\overline Q$ defined on 
 $\Z _+\times \mathscr {B} (\Z _+) $, we set 
 
 \[
\overline Q (y, B) = I_{\{0 \in B \}}, \ y \notin f(S).
\]

Then from the change of variables formula for the Lebesgue integral it
follows that
the last integral in \eqref{prance} allows the representation

$$
\int _S Q(x_{n-1} , dx_n) h_n(f(x_n)) = \
\int _{\Z _+} \overline Q (f(x_{n-1}) , dz_n) h_n(z_n).
$$

Likewise, we set $z_{n-1} = f(x_{n-1})$ in the next to last integral:

\[
\int _S Q(x_{n-2} , dx_{n-1}) h_n(f(x_{n-1}))
\int _S Q(x_{n-1} , dx_n) h_n(f(x_n)) = 
\]
\[ 
 \int _S Q(x_{n-2} , dx_{n-1}) h_n(f(x_{n-1}))
\int _{\Z _+} \overline Q (f(x_{n-1}) , dz_n) h_n(z_n) = 
\]
$$
\int _{\Z _+} \overline Q(f(x_{n-2}) , dz_{n-1}) h_n(z_{n-1})
\int _{\Z _+} \overline Q (z_{n-1} , dz_n) h_n(z_n).
$$
Further proceeding, we get
\[
\int _S Q(x_0 , dx_1) h_1(f(x_1)) \int _S Q(x_1 , dx_2) h_2(f(x_2))...
 \int _S Q(x_{n-1} , dx_n) h_n(f(x_n)) =
\]
\[
\int _{\Z _+} \overline Q(z_0 , dz_1) h_1(z_1) \int _{\Z _+} \overline Q(z_1 , dz_2) h_2(z_2)...
 \int _{\Z _+} \overline Q(z_{n-1} , dz_n) h_n(z_n),
\]
where
$z_0 = f (x _0 )$.

Thus,

\[
E_s \left[ \prod\limits _{j=1} ^k 
 h_j(Y_{n+j}) \mid \mathscr{F} ^ X _n  \right] =
 \]
 \[
\int _{\Z _+} \overline Q(f(X_0) , dz_1) h_1(z_1) 
\int _{\Z _+} \overline Q(z_1 , dz_2) h_2(z_2)...
 \int _{\Z _+} \overline Q(z_{n-1} , dz_n) h_n(z_n).
\]

  This equality and  \eqref{penultimate} imply that
  $Y$ is a Markov chain.

  \begin{rmk5}
   The kernel $\overline Q$ and the chain $f(X_n)$ are related:
   for all $s\in S$, $n,m \in \N$ and $M \subset \N$,
   
   \[
    P_s \{ f(X_{n+1}) \in M \mid f(X_{n}) = m \} = \overline Q(m,M)
   \]
whenever $P_s \{  f(X_{n+1}) = m \} >0$.
   Informally, one may say that 
   $\overline Q$ is the transition probability kernel for 
   the chain $\{f(X_n)\}_{n\in \Z _+}$.
  \end{rmk5}

  \begin{rmk2}  
   Clearly, this result holds for a Markov chain which is
   not necessarily defined on a canonical state space,
   because the property of a process to be a Markov chain
   depends on its distribution only.
  \end{rmk2}

\section{A birth-and-death process in the space of finite configurations: construction and
basic properties}

We would like to construct a Markov process 
in the space of finite configurations $\CYRG _0 (\R^{d})$, 
with a heuristic generator of the form

\begin{align} \label{the generator}
L F (\eta) = \int\limits _{x \in \R^d} b(x, \eta) [F(\eta \cup {x}) -F(\eta)] dx + 
\sum\limits _{x \in \eta} d(x, \eta ) (F(\eta \setminus {x}) - F(\eta)).
\end{align}
for
 $F$ in an appropriate domain. We call the functions  
  $b:\R ^d \times \Gamma _0 (\R^d) \rightarrow [0;\infty)$ and
 $d:\R ^d \times \Gamma _0 (\R^d) \rightarrow [0;\infty)$
 the \textit{birth rate coefficient} and
 the \textit{death rate coefficient}, respectively. 
 Theorem \ref{Core thm} summarizes the main results obtained in this section.

 To construct a spatial
 birth-and-death process, we consider the stochastic equation with Poisson noise

\begin{equation} \label{se}
\begin{split}
\eta _t (B) = \int\limits _{B \times (0;t] \times [0; \infty ] }
I _{ [0;b(x,\eta _{s-} )] } (u) dN_1(x,s,u) \\
 - \int\limits _{\Z \times (0;t] \times [0; \infty ) } I_{\{x_i \in \eta _{r-} \cap B \}}
I _{ [0;d(x_i,\eta _{r-} )] } (v)  d  N _2 (i,r,v) + \eta _0 (B),
\end{split}
\end{equation}
where
$(\eta _t)_{t \geq 0}$ is a suitable cadlag $\Go$-valued
stochastic process, the ``solution'' of the equation,
$B \in \mathscr{B} (\R ^d) $ is a Borel set, $N_1$
is a  Poisson point process on $ \R^d \times \R _+ \times \R _+ $
with intensity $dx \times ds \times du $,
$  N_2$ is a Poisson point process on 
$ \Z  \times \R _+ \times \R _+$ with intensity $\# \times dr \times dv$ ;
 $ \eta _0$
is a (random) finite initial configuration,
 ${b,d: \R ^d \times \Gamma _0 (R^d) \rightarrow [0;\infty)}$ are functions
measurable with respect to the product $\sigma$-algebra 
$\mathscr{B} (\R) \times \mathscr{B} ( \CYRG _0 (\R)) $, and
the sequence 
$\{...,x_{-1},x_{0},x_1,... \}$ is related to 
$(\eta _t)_{t \in [0;\infty]}$, as described in Section \ref{voranbringen}.
We require the
processes $N_1,  N _2, \eta _0$ to be independent of each other.
Equation \eqref{se}
is understood in the sense that the equality
holds a.s. for every bounded
$B \in \mathscr{B} (\R ^d) $ 
and $t \geq 0$.

As it was said in the preliminaries on Page \pageref{identify conf meas},
 we identify a finite configuration with a finite simple counting measure, so that
a configuration $\gamma$ acts as a measure in the following way:

\[
\gamma (A)  = |\gamma \cap A|, \ \ \ A \in \mathscr{B}(\R ^d).
\]
 
We will treat an element of 
 $\Go$ both as a set and as a counting measure, 
 as long as this does not lead to ambiguity. 
 An appearing of a new point
will be interpreted as a birth, 
and a disappearing will be interpreted as a death. We will refer to points
of $\eta _t$ as particles.

Some authors write $\tilde d(x,\eta \setminus x)$ where we write $d(x,\eta )$,
so that \eqref{the generator} translates to

\begin{align} 
L F (\eta) = \int\limits _{x \in \R^d} b(x, \eta) [F(\eta \cup {x}) -F(\eta)] dx + 
\sum\limits _{x \in \eta} \tilde d(x, \eta \setminus x ) (F(\eta \setminus {x}) - F(\eta)),
\end{align}
see
e.g. \cite{Preston}, \cite{Semigroupapproach}.

These settings are formally equivalent: the relation between $d$ and $\tilde d$ is 
given by
\[
 d(x,\eta) = \tilde d(x, \eta \setminus x), \ \ \ \eta \in \Go, x \in \eta,
\]
or, equivalently,
\[
 d(x,\xi \cup x) = \tilde d(x, \xi), \ \ \ \xi \in \Go, x \in \R ^d \setminus \xi.
\]

The settings used here appeared in \cite{HolleyStroock}, \cite{GarciaKurtz}, etc.

 We define the \emph{cumulative death rate} at $\zeta$ 
 by 
 
 \begin{equation}\label{cumulative death rate}
  D(\zeta) = \sum\limits _{x \in \zeta} d(x, \zeta),
 \end{equation}
 and
 the \emph{cumulative birth rate}
 by
  \begin{equation}\label{cumulative birth rate}
  B(\zeta) = \int\limits _{x \in \R^d} b(x, \zeta) dx.
 \end{equation}

\begin{def weak solution} \label{weak solution}
 A \emph{(weak) solution} of equation \eqref{se} is a triple 
 $(( \eta _t )_{t\geq 0} , N_1 ,  N _2)$, $(\Omega , \mathscr{F} , P) $,
 $(\{ \mathscr {F} _t  \} _ {t\geq 0}) $, where 

  (i) $(\Omega , \mathscr{F} , P)$ is a probability space, 
and $\{ \mathscr {F} _t  \} _ {t\geq 0}$ is an increasing, right-continuous
 and complete filtration of sub - $\sigma$ - algebras of $\mathscr {F}$,

  (ii)   $N_1$ is a
 Poisson point process on $\R ^d \times \R _+ \times \R _+$  with  
 intensity  $dx \times ds \times du $,

  (iii) $ N _2$ is a
 Poisson point process on $\Z  \times \R _+ \times \R _+$  with 
 intensity $\# \times ds \times du $,

  (iv) $ \eta _0  $ is a random  $\mathscr {F} _0$-measurable
  element in $\Go$,

  (v) the processes $N_1 ,  N _2$ and $\eta _0$ are independent,
 the processes $N_1$ and $ N _2$ are
  compatible with $\{ \mathscr {F} _t  \} _ {t\geq 0} $,

  (vi) $( \eta _t )_{t\geq 0} $ is a cadlag $\Go$-valued process
adapted to $\{ \mathscr {F} _t  \} _ {t\geq 0} $, $\eta _t \big| _{t=0} = \eta _0$,

  (vii) all integrals in \eqref{se} are well-defined, and
  
  (viii) equality \eqref{se} holds a.s. for all $t\in [0;\infty]$
  and all bounded Borel sets $B$, with 
  $\{ x_m\}_{m \in \Z}$ being the sequence related to $( \eta _t )_{t\geq 0} $.

\end{def weak solution}

Note that due to Statement \ref{distributions}
 item (viii) of this definition is a statement about the joint distribution of 
 $(\eta _t) , N_1 ,  N _2$.

 Let 
 \begin{align*}
   \mathscr{C} ^{0} _t =  \sigma \bigl\{ & \eta_0 , 
 N_1(B,[0;q],C),   N _2(i,[0;q],C);\\ &
B\in \mathscr{B} (\R^d), C\in \mathscr{B} (\R_+), q\in [0;t], i\in \Z \bigr\},
 \end{align*}
 and
 let $\mathscr{C} _t$ be the completion of $\mathscr{C} ^{0} _t$ under $P$.
 Note that $\{ \mathscr{C} _t \}_{t\geq 0} $ 
 is a right-continuous filtration, 
see Remark \ref{Missstand}.

\begin{def strong solution} \label{strong solution}
 A solution  of \eqref{se} is called \emph{strong}
 if $( \eta _t )_{t\geq 0} $ is adapted to 
$(\mathscr{C} _t, t\geq 0)$.

\end{def strong solution}

\begin{3rmk4}
 In the definition above we considered solutions as processes indexed by 
 $ t\ \in[0;\infty)$. The reformulations for 
the case  $ t \in [0;T]$, $0<T<\infty$, are straightforward. 
This remark applies to the results below, too.
\end{3rmk4}

Sometimes only the solution process (that is, $( \eta _t )_{t\geq 0} $)
 will be referred to as a (strong or weak) solution, when all the
 other structures are clear from the context.

We will say that the \textit{existence of strong solution} holds, 
if on any probability space with given
$N_1 , N _2, \eta _0$, satisfying (i)-(v) of Definition
\eqref{weak solution}, there exists a strong solution.

\begin{def pathwise uniqueness}
 We say that pathwise uniqueness holds for equation \eqref{se} and an
initial distribution $ \nu $ if, whenever the triples 
$(( \eta _t )_{t\geq 0} , N_1 ,  N _2)$, $(\Omega , \mathscr{F} , P) $,
 $(\{ \mathscr {F} _t  \} _ {t\geq 0}) $ and 
 $(( \bar \eta _t )_{t\geq 0} , N_1 ,  N _2)$, $(\Omega , \mathscr{F} , P) $,
 $(\{ \bar{\mathscr {F}} _t  \} _ {t\geq 0}) $ are weak solutions of \eqref{se} with 
$P \{ \eta _0 = \bar{\eta} _0 \} = 1 $ and $Law (\eta) = \nu $, 
we have $P \{ \eta _t = \bar{\eta} _t , t \in [0;T] \} = 1$
(that is, the processes $\eta , \bar{\eta} $ are indistinguishable).

\end{def pathwise uniqueness}

We assume that the birth rate $b$ satisfies 
the following conditions: sublinear growth on
the second variable in the sense 
that
\begin{equation} \label{sublinear growth for b} 
\int\limits _{\R ^d} {b}(x, \eta ) dx \leq c_1|\eta| +c_2,
\end{equation}
and let $d$ satisfy
\begin{equation} \label{condition on d}
\forall m \in \N : \sup\limits _{x\in \R^d,|\eta| \leq m} d(x, \eta) < \infty.
\end{equation}

 We also assume that

\begin{equation} \label{condition on eta _0}
 E |\eta _0| < \infty.
\end{equation}

By a non-random initial condition we understand an initial condition
with a distribution, concentrated at one point:
for some $\eta' \in \CYRG _0 (\R^d)$,
$P \{ \eta _0 = \eta'   \} =1$.

From now on, we work on some filtered probability space
$(\Omega , \mathscr{F}, (\{ \mathscr {F} _t  \} _ {t\geq 0}) , P) $.
On this probability space, 
the Poisson point processes $N_1 ,  N _2$ and $\eta _0$ are defined,
so that the whole set-up satisfies
(i)-(v) of Definition \ref{weak solution}.

Let us now consider the equation

\begin{equation} \label{pure birth}
\overline{\eta} _t (B) = \int\limits _{B \times (0;t] \times [0; \infty ] }
I _{ [0;\overline{b}(x,\overline{\eta} _s )] } dN(x,s,u) + \eta _0 (B),
\end{equation}
where
$\overline{b}(x, \eta ) := \sup\limits _{\xi \subset \eta} b(x, \xi) $.
Note that $\overline{b}$ satisfies sublinear growth condition
\eqref{sublinear growth for b}, if $b$ satisfies it.

This equation
is of the type \eqref{se} 
(with $\overline b$ being the birth rate coefficient,
and the zero function being the death rate coefficient), and all definitions 
of existence and uniqueness of solution
are applicable here. Later a unique solution
of \eqref{pure birth} will be
used as
 a majorant of a solution to \eqref{se}.

\begin{3lem2} \label{pure birth 3lem}
Under assumptions \eqref{sublinear growth for b} and \eqref{condition on eta _0},
strong existence and pathwise uniqueness hold for equation \eqref{pure birth}.
The unique solution $(\bar \eta _t)_{t\geq 0}$ satisfies 
\begin{equation}\label{beschwoeren}
 E|\bar \eta _t | < \infty, \ \ \  t \geq 0.
\end{equation}
\end{3lem2}

\textbf{Proof}. 
For $\omega \in  \{ \int\limits _{\R^d }
\overline{b}(x,\eta _0) dx =0 \} $, set $\zeta _t \equiv \eta _0$,
$\sigma _n = \infty$, $n \in \N$.

For $\omega \in F:= \{ \int\limits _{\R^d }
\overline{b}(x,\eta _0) dx >0 \}$,
 we define the sequence of random pairs $\{(\sigma _n, \zeta _{\sigma _n}) \}$, where
\[
\sigma _{n+1}= \inf \{ t>0 : \int\limits _{\R^d \times (\sigma _n; \sigma _n +t] \times [0;\infty)}
I_{[0; \overline b(x,\zeta _{\sigma _n})]} (u) dN_1(x,s,u) >0 \}+ \sigma _n, \ \ \sigma _0 = 0,
\]
and
\[
 \zeta _{0} = \eta _0, \ \ \ 
\zeta _{\sigma _{n+1}} = \zeta _{\sigma _n} \cup \{ z_{n+1} \}
 \]
for $z_{n+1} = \{x\in \R^d: N_1 (x,\sigma _{n+1}, [0; \overline b(x,\zeta _{\sigma _n})]) >0  \}$.
From \eqref{travesty} it follows that
the points 
$z_n$ are uniquely determined almost surely on $F$. 
Moreover, 
 $\sigma _{n+1} > \sigma _n$ a.s., and $\sigma _n$
 are finite a.s. on $F$
 (particularly because $\overline{b}(x,\zeta _{\sigma _n}) \geq \overline{b}(x,\eta _0)$).
For $\omega \in F$, we define $\zeta _t = \zeta _{\sigma _n}$
for $t\in [\sigma _n; \sigma _{n+1})$. Then by induction on $n$ 
it follows that $\sigma _n $ is a stopping time for each $n \in \N$, and 
$\zeta _{\sigma _n}$ is $\mathscr{F} _{\sigma _n} \cap F $-measurable. 
By direct substitution we see that
 $(\zeta _t)_{t\geq 0}$ is a strong
 solution for \eqref{pure birth} on the time interval
$t\in [0; \lim\limits _{n\to \infty} \sigma _n)$. 
Although we have not defined what is a solution,
or a strong solution, on a 
random time interval, we do not discuss it here.
Instead we are going to show that 
\begin{equation}\label{staunch}
\lim\limits _{n\to \infty} \sigma _n = \infty \ \ \ \textrm{a.s.}
\end{equation}
This relation is evidently true on the complement of $F$.
If $P(F)=0$, then \eqref{staunch} is proven.

If $P(F)>0$, define a probability measure on $F$, 
$Q(A) = \frac{P(A)}{P(F)}$, $A \in \mathscr{S} := \mathscr{F} \cap F$,
and define $ \mathscr{S}_t =\mathscr{F}_t \cap F $.

The process $N_1$ is independent of $F$, therefore it is a
Poisson point process on $(F,\mathscr{S}, Q)$ with the same intensity,
compatible with $\{ \mathscr{S}_t \}_{t\geq 0}$.
From now on and until other is specified, we work on the filtered
probability space $(F,\mathscr{S},\{ \mathscr{S}_t \}_{t\geq 0}, Q)$.
We use the same symbols for random processes and random variables, 
having in mind that 
we consider their restrictions to $F$.

The process $(\zeta _t)_{t\in [0; \lim\limits _{n\to \infty} \sigma _n)}$
has the Markov property, because 
the process $N_1$ has the strong Markov property and independent increments. 
Indeed, conditioning on $ \mathscr{S}_{\sigma_n}$,
\[
 E \bigl[  I_{\{\zeta _{\sigma _{n+1}} = \zeta _{\sigma _{n}} \cup x 
 \text{ for some } x \in B \}} \mid \mathscr{S}_{\sigma_n}\bigr]=
 \frac{\int\limits _{ B} \overline b(x, \zeta _{\sigma _{n}}) dx}{\int\limits _{\R^d }
\overline{b}(x,\zeta _{\sigma _n}) dx },
\]
thus the chain $\{ \zeta _{\sigma _{n}}\}_{n \in Z_+}$ is a Markov chain,
and, given $\{ \zeta _{\sigma _{n}}\}_{n \in Z_+}$,
$\sigma _{n+1} - \sigma _n $ are distributed exponentially:
\[
E\{ I_{\{\sigma _{n+1} - \sigma _n >a\}}  \mid \{ \zeta _{\sigma _{n}}\}_{n \in Z_+}\} 
= \exp \{ - a \int\limits _{\R^d }
\overline{b}(x,\zeta _{\sigma _n}) dx \}.
\]
Therefore, the random variables $\gamma _n = (\sigma _{n} - \sigma _{n-1}){(\int\limits _{\R^d }
\overline{b}(x,\zeta _{\sigma _n}) dx)}$ 
constitute a 
sequence of independent 
 random variables exponentially distributed
with parameter $1$, independent 
of $\{ \zeta _{\sigma _{n}}\}_{n \in Z_+}$.
 Theorem 12.18 in \cite{KallenbergFound} (see Page \pageref{synthesis}
of this article) implies  that
$(\zeta _t)_{t\in [0; \lim\limits _{n\to \infty} \sigma _n)}$
is a pure jump type Markov process.

The jump rate of $(\zeta _t)_{t\in [0; \lim\limits _{n\to \infty} \sigma _n)}$
is given by
\[
 c(\alpha) = \int\limits _{\R^d }
\overline{b}(x,\alpha) dx .
\]
Condition \eqref{sublinear growth for b} implies that 
$c(\alpha) \leq c_1 |\alpha| + c_2$. Consequently,
\[
c (\zeta _{\sigma _n} ) \leq c_1 |\zeta _{\sigma _n}| + c_2 = 
c_1 |\zeta _0| +c_1 n  + c_2.
\]

We see that $\sum _n \frac{1}{c (\zeta _{\sigma _n} )} = \infty$ a.s., 
hence Proposition 12.19 in \cite{KallenbergFound} (given in 
Section
\ref{Pure jump type MP})
implies that $\sigma _n \to \infty$.

Now, we return again to our initial probability space 
$(\Omega,\mathscr{F},\{ \mathscr{F}_t \}_{t\geq 0}, P)$.

Thus, we have existence of a strong solution. 
Uniqueness follows by induction on jumps of the process. 
Indeed, let $( \tilde{\zeta} _t )_{t\geq 0} $ be another solution of \eqref{pure birth}. From (viii)
 of Definition \ref{weak solution} and equality

$$
\int\limits _{\R^d \times (0;\sigma _1) \times [0; \infty ] }
I _{ [0;\overline{b}(x, {\eta} _0 )] } dN_1(x,s,u) = 0 ,
$$
one
can see that $P \{ \tilde{\zeta} \text{\ has a birth before \ }
\sigma _1  \} = 0$. At the same time, equality 

$$
\int\limits _{\R^d \times \{ \sigma _1 \} \times [0; \infty ] }
I _{ [0;\overline{b}(x, {\eta} _0 )] } dN_1(x,s,u) = 1 ,
$$
which holds a.s., 
yields that $\tilde{\zeta}$ has a birth at the moment $\sigma _1$, and in the same point 
of space at that. Therefore, $\tilde{\zeta}$ coincides with $\zeta$ up to $\sigma _1$ a.s. Similar reasoning shows
that they coincide up to $\sigma _n$ a.s., and, because $\sigma _n \to \infty$ a.s., 

$$P \{ \tilde{\zeta} _t = {\zeta} _t \text{\ for all \ } t\geq 0 \} = 1$$

 Thus, pathwise uniqueness holds. The constructed
 solution is strong.
 
Now we turn our attention to \eqref{beschwoeren}.
We can write 

\begin{gather}\label{gibberish}
|{\zeta} _t| = 
 |{\eta} _0|+  \sum\limits _{n=1}^{\infty} I \{ |\zeta _t| - |\eta _0| \geq n \} \notag
 \\  = 
  |{\eta} _0|+  \sum\limits _{n=1}^{\infty} I \{ \sigma _n \leq t \}.
\end{gather}

Since $\sigma _n = \sum\limits _{i=1}^n \frac{\gamma _i}{\int\limits _{\R^d }
\overline{b}(x,\zeta _{\sigma _i}) dx}$, we have

\[
 \{ \sigma _n \leq t \} = \{ \sum\limits _{i=1}^n \frac{\gamma _i}{\int\limits _{\R^d }
\overline{b}(x,\zeta _{\sigma _i}) dx} \leq t \} \subset
\{ \sum\limits _{i=1}^n \frac{\gamma _i}{c_1|\zeta _{\sigma _i}| + c_2} \leq t \}
\]
\[
\subset \{ \sum\limits _{i=1}^n \frac{\gamma _i}{(c_1+c_2)(|\eta _0|+i)} \leq t \} =
\{Z _t -Z_0 \geq n \},
\]
where
$(Z_t)$ is the Yule process (see Page \pageref{Yule page})
with birth rate defined as follows:
$Z_t - Z _0 = n$ when 
\[
\sum\limits _{i=1}^n \frac{\gamma _i}{(c_1+c_2)(|\eta _0|+i)} \leq t 
< \sum\limits _{i=1}^{n+1} \frac{\gamma _i}{(c_1+c_2)(|\eta _0|+i)},
\]
and 
$Z_0 = |\eta _0|$. Thus, we have $| \zeta _t| \leq Z _t$ a.s.,
hence $E| \zeta _t| \leq  E Z _t < \infty$.
 $\Box$

\begin{3thm}\label{ex un} 
Under assumptions \eqref{sublinear growth for b}-\eqref{condition on eta _0},
 pathwise uniqueness and strong existence hold for equation \eqref{se}. 
 The unique solution $(\eta _t)$ is a pure jump type process satisfying 
 \begin{equation}\label{unfathomable}
 E| \eta _t | < \infty, \ \ \  t \geq 0.
\end{equation}
\end{3thm}

\textbf{Proof}. Let us define stopping times with respect to
 ${\{ \mathscr{F} _t , t\geq 0 \}} $,
$0= \theta _0 \leq \theta _1 \leq \theta _2 \leq \theta _3 \leq ...$,
and the sequence of (random) configurations 
$\{ \eta _{\theta _j} \} _{j \in \N}$
as follows: as long as 
\[
 B({\eta} _{\theta _{n}})+D({\eta} _{\theta _{n}})>0,
\]
we set
\[ 
\theta _{n+1} = \theta ^b _{n+1} \wedge \theta ^d _{n+1} + \theta _n,
\]
\[
\theta ^b _{n+1}= \inf \{ t>0 : \int\limits _{\R^d \times (\theta _n; \theta _n +t] \times [0;\infty)}
I_{[0; b(x,\eta _{\theta _n})]} (u) dN_1(x,s,u) >0 \},
\]
\[
\theta ^d _{n+1}= \inf \{ t>0 : \int\limits _{(\theta _n; \theta _n +t]
 \times [0; \infty ) } I _{\{x_i\in  \eta _{\theta _n}\}}
I _{ [0;d(x_i,\eta _{\theta _n} )] } (v) dN_2 (i,r,v) >0 \},
\]
${\eta} _{\theta _{n+1}}
=  {\eta} _{\theta _n} \cup \{ z_{n+1} \}$ if $\theta ^b _{n+1} \leq \theta ^d _{n+1}$,
where
$\{z_{n+1}\} = \{z\in \R^d: N_1 (z,\theta _n + \theta ^b _{n+1}, \R_+) >0  \}$; 
${\eta} _{\theta _{n+1}} =  {\eta} _{\theta _n} \setminus \{ z_{n+1} \}$
if  $\theta ^b _{n+1} > \theta ^d _{n+1}$, where
$\{z_{n+1}\} = \{x_i \in {\eta} _{\theta _n}: N_2  (i,\theta _n + \theta ^d _{n+1}, \R_+) >0  \}$;
the configuration $\eta _{\theta _0} = \eta _0$
is the initial condition of \eqref{se}, ${\eta} _t = {\eta} _{\theta _n} $ for 
$t \in [\theta _n ; \theta _{n+1} )$, $\{x_i \}$
is the sequence related to $({\eta} _t)_{t \geq 0}$.
 Note that 
\[
P\{  \theta ^b _{n+1} = \theta ^d _{n+1} \text{\ for some \ } n \mid
  B({\eta} _{\theta _{n}})+D({\eta} _{\theta _{n}}) >0 \} = 0,
\]
 the points $z_n$ are a.s. uniquely determined, and 
 \[
 P\{  z_{n+1} \in {\eta} _{\theta _n}\mid \theta ^b _{n+1} \leq \theta ^d _{n+1} \} =0.
 \]
If for some $n$
\[
 B({\eta} _{\theta _{n}})+D({\eta} _{\theta _{n}})=0,
\]
then we set $\theta _{n+k} =\infty$, $k\in \N$, and
$\eta _{t} =  {\eta} _{\theta _{n}}$, $t \geq {\theta _{n}}$.

As in the
proof of Proposition \ref{pure birth 3lem}, $({\eta} _t)$ is a strong solution of \eqref{se},
$t \in [0;\lim_n \theta _n)$.

Random variables $\theta _n, n \in \N$, are stopping times with respect to
the filtration ${\{ \mathscr{F} _t , t\geq 0 \}} $. 
Using the strong
Markov property of a Poisson point process, we see that,
on $\{ \theta _n < \infty \} $, the conditional distribution of $\theta ^b _{n+1}$ given
$\mathscr {F} _{\theta _n}$ is $\exp (\int\limits _{\R^d} b(x,\eta _{\theta _n}) dx)$,
and the conditional distribution of $\theta ^d _{n+1}$ given
$\mathscr {F} _{\theta _n}$ is $\exp ( \sum\limits _{x \in \eta _{\theta _n}} d(x,\eta _{\theta _n} )) $. 
In particular,
$\theta ^b _n, \theta ^d_n > 0$, $n \in \N$, and the process $(\eta _t)$ is 
of pure jump type. 

Similarly to the proof of Proposition \ref{pure birth 3lem}, 
one can show by induction on $n$ that equation \eqref{se} 
has a unique solution on $[0; \theta _n]$. Namely, each two
solutions coincide on $[0; \theta _n]$ a.s.
Thus, any solution coincides with $({\eta} _t)$ a.s. for all 
$t  \in [0; \theta _n] $.

 Now we will show that $\theta _n \to \infty$ a.s. as $n \to \infty$.
Denote by $\theta _k '$ the moment of the $k$-th birth. It is sufficient to show that
$\theta _k ' \to \infty$, $k\to \infty$, 
because only finitely many deaths may occur between any two births, since there
are only finitely particles. By induction
on $k '$ one may see that 
$\{ \theta _k' \} _{k' \in \N} \subset \{ \sigma _i \} _{i \in \N}$, where
$\sigma _i$ are the moments of births of $(\overline{\eta} _t)_{t\geq 0}$,
the solution of \eqref{pure birth},
 and $\eta _t \subset \overline{\eta} _t$
 for all $t \in [0;\lim_n \theta _n)$.
For instance, let us show that $(\overline{\eta} _t)_{t\geq 0}$ has a birth at
 $\theta _1 '  $. We have  $\overline{\eta} _{\theta _1 ' -} \supset \overline{\eta} _{0} = \eta _0$, and
$\eta _{\theta _1 ' -} \subset \eta _{t} \mid _{t=0} = \eta _0 $, 
hence for all $x \in \R^d$

$$\overline{b} (x, \overline{\eta} _{\theta _1 ' -} ) \geq \overline{b} (x, \eta _{\theta _1 ' -} ) \geq b (x, \eta _{\theta _1 ' -} )$$

The latter implies that at time moment $\theta _1 '  $ a birth
 occurs for the process $(\overline{\eta} _t)_{t\geq 0}$ in the same point.
 Hence, $\eta _{\theta _1 '} \subset \overline{\eta} _{\theta _1 ' } $,
 and we can go on.
Since $\sigma _k \to \infty$ as
$k \to \infty$, we also have $\theta _k ' \to \infty$, and therefore $ \theta _n \to \infty$, $n\to \infty$. 

 Since $\eta _t \subset \overline{\eta} _t$ a.s., Proposition \ref{pure birth 3lem} 
 implies \eqref{unfathomable}. 
 $\Box$

 In particular, for any time $t$ the integral
\[
\int\limits _{\R ^d \times (0;t] \times [0; \infty ] }
I _{ [0;b(x,\eta _{s-} )] } (u) dN_1(x,s,u)
\]
 is finite a.s.

\begin{3rmk3}\label{encompass}
Let $\eta _0$ be a non-random initial condition, $\eta _0 \equiv \alpha$,
$\alpha \in \Go$.
The solution of \eqref{se} with $\eta _0 \equiv \alpha$
will be denoted as $(\eta(\alpha, t))_{t\geq 0}$.
Let $P_\alpha$ be
the push-forward of $P$ under the mapping

\begin{equation}
  \Omega \ni \omega \mapsto (\eta(\alpha, \cdot)) \in D_{\Go}[0;T].
\end{equation}

 From the proof one may derive that, for fixed $\omega \in \Omega$, 
constructed unique solution is jointly measurable in $(t, \alpha)$.
Thus, the family $\{ P_\alpha \}$
 of probability measures on $D_{\Go}[0;T]$
is measurable in $\alpha$.
 We will often use formulations
related to the probability space
$(D_{\Go}[0;T],\mathscr{B}(D_{\Go}[0;T]), P_{\alpha})$;
in this case, coordinate mappings will be denoted 
by $\eta_t$,

\[
 \eta _t (x) = x(t), \ \ \   x \in D_{\Go}[0;T].
\]

The processes $(\eta _t) _{t \in  [0;T]}$
and $(\eta(\alpha, \cdot))_{t \in  [0;T]} $
have the same law (under $P_{\alpha} $ and $P$, respectively).
As one would expect, the family
of measures $\{P_{\alpha}, \alpha \in \Go \}$
is a Markov process, or a Markov family of probability measures;
see Theorem \ref{Markov property}
below. For a measure $\mu$ on $\Go$,
 we define
 
 \[
  P_{\mu} = \int P_{\alpha} \mu (d \alpha).
 \]

We denote by $E_{\mu}$ the expectation under $P_{\mu}$.

\end{3rmk3}

\begin{3rmk2}\label{3rmk2}
Let $b_1,d_1$ be another pair of birth and death coefficients, 
satisfying all conditions imposed on 
$b$ and $d$. 
Consider a unique solution $(\tilde{\eta} _t)$ of \eqref{se} with coefficients 
$b_1, d_1$ instead of
$b,d$, but with the same initial condition $\eta _0$ and all the
other underlying structures. \textit{If for all }
$\zeta \in D$, \textit{where}
$D \in \mathscr{B}( \CYRG _0 (\R^d))$ , 
$b_1 (\cdot ,\zeta) \equiv b (\cdot ,\zeta) $, $d_1 (\cdot ,\zeta) \equiv d (\cdot ,\zeta) $,
 \textit{then} 
$\tilde{\eta} _t = \eta _t$ \textit{\ for all } $t\leq \inf \{s \geq 0: \eta_s \notin D\} = 
\inf \{s \geq 0: \tilde{\eta}_s \notin D\} $.
This may be proven in the same way as the theorem above.

\end{3rmk2}

\begin{3rmk8}
Assume that all the conditions
of Theorem \ref{ex un} are fulfilled 
except Condition  \eqref{condition on eta _0}.
Then we could not claim that 
 \eqref{unfathomable} holds.
However, other conclusions of the theorem would hold. 
We are mostly interested in the case of a non-random 
initial condition, therefore we do not discuss 
the case when \eqref{unfathomable} is not satisfied.
\end{3rmk8}

\begin{3rmk9}\label{Functional dependence}
 We solved equation \eqref{se} 
 $\omega$-wisely. As a consequence,
 there is a functional 
 dependence of the solution
 process and the ``input'':
 the process $(\eta _t)_{t\geq 0}$
 is some function of 
 $\eta _0$, $N_1$ and $N_2$.
 Note that $\theta _n$ and $z_n$ from the proof 
 of Theorem \ref{ex un} are measurable 
 functions of $\eta _0$, $N_1$ and $N_2$
 in the sense that, e.g., 
 $\theta_1 = F _1 (\eta _0, N_1, N_2)$ a.s.
 for a measurable $F_1: \Go \times \CYRG(\R ^d \times \R _+ \times \R _+)
 \times \CYRG(\Z ^d \times \R _+ \times \R _+) \to \R _+$.
\end{3rmk9}

\begin{3lem1} \label{expectation}  
If $(\eta _t)_{t\geq 0}$ is a solution to equation \eqref{se},
 then the inequality
\[
E|\eta _t|< (c_2 t + E |\eta _0|)e^{c_1 t}
\]
holds for all $t>0$.

\end{3lem1}

Proof. We already know that
$E|\eta _t|$ is finite. Since $\eta _t $ satisfies equation \eqref{se} we have

\[
\eta _t (B) = \int\limits _{B \times (0;t] \times [0; \infty ] }
I _{ [0;b(x,\eta _{s-} )] } (u) dN_1(x,s,u) 
\]
\[
- \int\limits _{\Z \times (0;t] \times [0; \infty ) } I_{\{x_i \in \eta _{r-} \cap B \}}
I _{ [0;d(x_i,\eta _{r-} )] } (v)  d  N _2 (i,r,v) \leq
\] 
\[
 \int\limits _{B \times (0;t] \times [0; \infty ] }
I _{ [0;b(x,\eta _{s-} )] } (u) dN_1(x,s,u) + \eta _0 (B).
\]
For $B=\R^d$, 
taking expectation in the last inequality,
we obtain
\[
E |\eta _t| = E\eta _t (\R^d)  \leq E \int\limits _{\R^d \times (0;t] \times [0; \infty ] } I _{ [0;b(x,\eta _{s-} )] } (u) dN_1(x,s,u) +
 E \eta _0 (\R^d) =
\] 
\[
= E \int\limits _{\R^d \times (0;t] \times [0; \infty ] } I _{ [0;b(x,\eta _{s-} ) ]} (u) dxdsdu + E \eta _0 (\R^d) =
 E \int\limits _{\R^d \times (0;t] }  b(x,\eta _{s-} ) dxds + E \eta _0 (\R^d).
\]

Since $\eta$ is a solution of \eqref{se}, we have for all $s\in [0;t]$ almost surely $\eta _{s-} = \eta _s$. Consequently,
$E|\eta _{s-} |= E| \eta _s |$. Applying this and \eqref{sublinear growth for b}, we see that 

$$
E\eta _t (\R^d) \leq E \int\limits _{(0;t] } ( c_1 |\eta _{s-}| + c_2 ) ds + E \eta _0 (\R^d)
= c_1  \int\limits _{(0;t] }  E |\eta _s| ds + c_2 t + E \eta _0 (\R^d),
$$
so the
statement of the lemma follows from \eqref{condition on eta _0} 
and Gronwall's inequality. $\Box$

\begin{def joint uniqueness in law} \label{joint uniqueness in law}
 We say that \textit{joint uniqueness in law} holds for equation \eqref{se} with an initial
distribution $\nu$ if any two (weak) solutions $((\eta_t) , N_1 , N_2)$ and 
$((\eta_t)  ^{ \prime } , N_1  ^{\prime}, N _2  ^{\prime}  )$ of \eqref{se},
$Law(\eta _0)= Law(  (\eta _0)  ^{\prime})=\nu$, have the same joint distribution:

$$Law ((\eta_t) , N_1 , N_2)
= Law ((\eta_t)  ^{\prime} , N_1  ^{\prime}, N_2  ^{\prime}  ).$$

\end{def joint uniqueness in law}

The following corollary is a consequence of Theorem 
\ref{ex un} and Remark \ref{Functional dependence} .

\begin{3cor1}
Joint uniqueness in law holds for equation \eqref{se} with initial
distribution $\nu$ satisfying

\[
 \int _{\Go} |\gamma| \nu (d \gamma) < \infty.
\]

\end{3cor1}

\begin{3rmk10}\label{ordering}
 We note here that altering the order of 
 the initial configuration does not change
 the law of the solution. We could replace 
 the lexicographical order with any other.
 To see this, note that if
 $\varsigma$ is a permutation of $\Z$ (that is,
 $\varsigma: \Z \to \Z$ is a bijection), then 
 the process $\tilde N_2$ defined by
\begin{equation}\label{permutation}
  \tilde N_2 (K, R,V) = N_2 (\varsigma K, R,V), \
  \ \ K \subset \Z, R,V \in \mathscr{B}(\R_+),
\end{equation}
has the same law as $N_2$, and is adapted 
to $\{\mathscr{F}_t \}_{t\geq 0}$, too. 
Therefore, solutions of \eqref{se} and of \eqref{se}
with $N_2$ being 
replaced by $\tilde N_2$ have the same law. But
replacing $N_2$ with $\tilde N_2$ in equation \eqref{se}
is equivalent to replacing
$\{ x_{-|\eta _0|+1},...,x_0,x_1,... \}$
with \\
$\{ x_{\varsigma^{-1}(-|\eta _0|+1)},...,
 x_{\varsigma ^{-1}(0)},x_{\varsigma ^{-1}(1)},... \}$.

\end{3rmk10}
Let $\nu$ be a distribution on $\Go$, and let $T>0$.
Denote by $\mathscr{L}(\nu, b, d, T)$ the law of the
restriction $(\eta _t)_{t \in [0;T]}$ of the unique solution 
$(\eta _t)_{t \geq 0}$
to 
\eqref{se} with an initial condition distributed according
to $\nu$. Note that $\mathscr{L}(\nu, b, d, T)$
is a distribution on $D_{\Go}([0;T])$.
As usually, the Markov property of a solution follows from uniqueness.

\begin{3thm3}\label{Markov property}
 The unique solution $(\eta _t) _{t\in [0;T]}$ of \eqref{se} is a Markov process.
\end{3thm3}

\textbf{Proof}. Take arbitrary $t'<t$, $t',t \in [0;T]$. Consider 
the equation

\begin{equation} \label{se t'<t}
\begin{split}
\xi _t (B) = \int\limits _{B \times (t';t] \times [0; \infty ] }
I _{ [0;b(x,\xi _{s-} )] } (u) dN_1(x,s,u) \\
 - \int\limits _{\Z \times (t';t] \times [0; \infty ) } I_{\{x_i ' \in \xi _{r-} \cap B \}}
I _{ [0;d(x_i ' ,\xi _{r-} )] }  d N _2 (i,r,v) + \eta _{t'}(B),
\end{split}
\end{equation}
where
the sequence $\{x_i ' \}$ is related to 
the process $(\xi _s)_{s \in [0;t]}$,
$\xi _s = \eta _s$.
The unique
solution of \eqref{se t'<t}
is $(\eta _s) _{s\in [t';t]}$.
As in the proof of Theorem \ref{ex un}
we can see that $(\eta _s) _{s\in [t';t]}$
is measurable with respect to the filtration 
generated by the random variables 
 $N_1(B, [s;q], U)$,
$ N_2 (i, [s;q],U)$, 
and $\eta _{t'}(B)$, where
$B \in \mathscr{B} (\R ^d)$, $i \in \Z$,
$t'\leq s \leq q \leq t$, $U \in \mathscr{B} (\R _+)$.
Poisson point process have independent increments,
hence
\[
 P \{ (\eta _t) _{t\in [s;T]} \in U \mid \mathscr{F}_s \} = 
 P \{ (\eta _t) _{t\in [s;T]} \in U \mid \eta _s \}
\]
almost surely. Furthermore, using
arguments similar to those in
Remark \ref{ordering}, we can conclude that 
$(\eta _s) _{s\in [t';t]}$ is distributed
according to $\mathscr{L}(\nu _{t'},b,d,t-t')$, where
$\nu _{t'}$ is the distribution of $\eta _{t'}$.
$\Box$

The following theorem sums up the results 
we have obtained so far.

\begin{3thm6}\label{Core thm}
 Under assumptions \eqref{sublinear growth for b}, \eqref{condition on d},
 \eqref{condition on eta _0}, equation \eqref{se} has a unique solution.
 This solution is a pure jump type Markov process.  The family of
 push-forward measures $\{P_{\alpha}, \alpha \in \Go \}$ defined
 in Remark \ref{encompass}
 forms a Markov process, or a Markov family of probability measures,
 on $D_{\Go}[0;\infty)$.
\end{3thm6}
\textbf{Proof}. The statement is a 
consequence of Theorem \ref{ex un}, Remark \ref{encompass} and
Theorem \ref{Markov property}. In particular, 
the Markov property of $\{P_{\alpha}, \alpha \in \Go \}$
follows from the statement given in the last sentence of 
the proof of Theorem \ref{Markov property}.
$\Box$

We call the unique solution of \eqref{se} (or, sometimes, the corresponding 
family of measures on $D_{\Go}[0;\infty)$) a
\textit{(spatial) birth-and-death Markov process}.

\begin{3rmk6}\label{domain of d}

We note that $d$
does not need to be defined on the whole space $\R ^d \times \Gamma _0 (\R^d)$.
The equation makes sense even if $d(x, \eta)$ is
defined on $\{ (x, \eta)\mid x \in \eta \}$. Of course,
any such function may be extended to a function
on $\R ^d \times \Gamma _0 (\R^d)$.
\end{3rmk6}

\subsection{Continuous dependence on initial conditions}

In order to prove the continuity of the distribution of the solution of \eqref{se}
 with respect to initial conditions, we make 
the following continuity assumptions on $b$ and $d$. 

\begin{3Assumpt}\label{Cont assumptions}
 Let $b,d$ be continuous with respect to both arguments. 
Furthermore, let the map 
\[
  \Go \ni \eta \mapsto b(\cdot , \eta) \in L^1(\R ^d)
\]
be continuous.

\end{3Assumpt}

In light of Remark \ref{domain of d},
let us explain what we understand by continuity of $d$ when
$d(x,\eta)$ is defined only on $\{ (x, \eta)\mid x \in \eta \}$. We
require that, whenever $\eta _n \to \eta$
and $\eta _n \ni z_n \to x \in \eta$, we also have
$d(z_n ,\eta _n) \to d (x, \eta )$.
Similar condition appeared in 
\cite[Theorem 3.1]{HolleyStroock}.

\begin{3thm2} \label{continuity of solution} 
 
Let the birth and death coefficients $b$ and $d$ satisfy the above continuity assumptions
\ref{Cont assumptions}. Then for every $T>0$
the map 

\[\CYRG _0 (\R^d) \ni \alpha  \mapsto 
Law \{ \eta(\alpha, \cdot) _\cdotp, \cdotp \in (0;T] \},
\]
which
assigns to a non-random initial condition $\eta _0 = \alpha $ the law of the
solution of equation \eqref{se} stopped at time $T$, is continuous.
\end{3thm2}
 \textbf{Remark}.
We mean continuity in in the space of measures on 
$D_{\Go} [0;T]$; see Page \pageref{D space}.

\textbf{Proof}. Denote by $\eta (\alpha, \cdotp )$ the solution of \eqref{se}, 
started from  $\alpha $. 
Let $\alpha _n \to \alpha$, $\alpha _n, \alpha \in \Go$,
$\alpha = \{x_0, x_{-1},...,x_{-|\alpha|+1} \}$,
$x_{0} \preccurlyeq x_{-1} \preccurlyeq  ... \preccurlyeq x_{-|\alpha|+1}$.
With no loss in generality we assume that $|\alpha _n|=|\alpha|$, 
$n \in \N$.
By Lemma \ref{caveat} we can label elements of $\alpha _n$,
$\alpha _n = \{x_0 ^{(n)}, x_{-1}^{(n)},...,x_{-|\alpha|+1}^{(n)} \}$,
so that $x_{-i}^{(n)} \to x_{-i}$, $i=0,...,|\alpha|-1$.
Taking into account Remark \ref{ordering}, we can assume
\begin{equation}\label{svelte}
 x_{0}^{(n)} \preccurlyeq x_{-1}^{(n)} \preccurlyeq  ... \preccurlyeq x_{-|\alpha|+1}^{(n)}
\end{equation}
 without loss of generality (in the sense that we do not have
 to use lexicographical order; not in the sense that we can make 
 $x_{0}^{(n)},x_{-1}^{(n)},...$ satisfy \eqref{svelte}
 with the lexicographical order).

 We will show that 

 \begin{equation}\label{accent}
  \sup\limits _{t\in [0;T] }  dist( \eta (\alpha , t ), \eta ( \alpha _n , t ) )
 \overset{p}{\to} 0, \ \ n \to \infty .
  \end{equation}

Let $ \{ \theta _i \} _{i \in \N} $ be the moments of jumps of process $\eta (\alpha , \cdotp )$. 
Without loss of generality, 
assume that $d(x,\alpha)> 0$, $x \in \alpha$, and 
$||b(\cdot , \alpha)||_{L^1}>0$, $L^1 := L^1(\R ^d)$ 
(if some of these inequalities are not fulfilled,
the following reasonings should be changed insignificantly).

Depending on whether 
a birth or a death occurs at $\theta _1$, we have  either 
\begin{equation}\label{that is to say}
  N_1 (\{x_1\}\times \{ \theta _1 \}\times [0;b(x_1, \eta_0)]) = 1
\end{equation}
or for some $x_{-k} \in \alpha$
\[
 N_2 (\{-k\}\times\{ \theta _1 \}\times [0;d(x_{-k}, \alpha)]) = 1.
\]
The probability
of last two equalities holding simultaneously
is zero, hence we can neglect this event. 
In both cases $N_1 (x_1,\{ \theta _1 \}, \{ b(x_1, \alpha)\}) = 0$, 
$N_2 (-k,\{\theta _1 \}, \{ d(x_{-k}, \alpha)\}) =0$ a.s. 
We also have
\[
 N_1 (\R ^d \times [0;\theta _1)\times [0;b(x, \alpha)]) = 0,
\]
and 
for all $j \in {0,1,...,|\alpha| -1}$
\[
N_2 (\{-j\} \times [0;\theta _1) \times [0;d(x_{-j}, \alpha)]) = 0.
\]

Denote 
\[
 m:= b(x_1,\alpha) \wedge \min \{d(x,\alpha): x \in \alpha \} \wedge
||b(\cdot , \alpha)||_{L^1} \wedge 1
\]
and
fix $\varepsilon >0$.
Let $\delta _1 >0$ be so small that for 
$\nu \in \Go$, $\nu = \{x'_0, x'_{-1},...,x'_{-|\alpha|+1} \}$,
$|x_{-j} - x'_{-j}|\leq \delta _1$
the inequalities
\[
|d( x'_{-j},\nu) - d( x_{-j}, \alpha)| < \varepsilon m , 
\]
and
\[
||b(\cdot ,\nu) - b(\cdot , \alpha)||_{L^1} < \varepsilon m
\]
hold. 
Then we may estimate

\begin{equation}\label{squash}
P \big\{ \int\limits _{\R^d \times [0;\theta _1 ) \times [0; \infty ] }
I _{ [0; b(x,\nu )] } (u) dN_1(x,s,u) \geq 1 \big\} < \varepsilon.
\end{equation}
and
\begin{equation}\label{proliferation}
P \big\{ \int\limits _{\Z \times [0;\theta _1 ) \times [0; \infty ] } I _{\{x'_{-i} \in \nu\}}
I _{ [0; d(x'_{-i},\nu )] } (v) dN_2 (i,r,v) \geq 1 \big\} < \varepsilon |\alpha|.
\end{equation}

Indeed, the random variable 

\begin{equation}
\tilde \theta := \inf\limits _{t >0} 
\{ \int\limits _{\R^d \times [0;t ) \times [0; \infty ] }
I _{ [0;0 \vee \{ b(x,\nu ) - b(x,\alpha\} )] } (u) dN_1(x,s,u) \geq 1 \} 
\end{equation}
is 
exponentially distributed with parameter 
$||(b(\cdot,\nu ) - b(\cdot,\alpha)) _+ ||_{L^1} < \varepsilon ||b(\cdot,\alpha) ||_{L^1}$.
By Lemma 
\ref{preclude}, 
\begin{equation}
P \{ \tilde \theta < \theta _1 \} < \frac{\varepsilon ||b(\cdot,\alpha) ||_{L^1}}
{||b(\cdot,\alpha) ||_{L^1}} = \varepsilon,
\end{equation}
which 
is exactly  \eqref{squash}.
Likewise, \eqref{proliferation} follows.

Similarly, the probability that the
same event as for 
$\eta (\alpha , \cdot)$ occurs at time $\theta _1$ 
for $ \eta (\nu , \cdot)$  is high.
Indeed, assume, for example,
that a birth occurs at
$\theta _1$, that is to say that  
\eqref{that is to say} holds. 
Once more using Lemma \ref{preclude}
we get

\[
 P \{ N_1 (\{x_1\}\times \{ \theta _1 \}\times [0;b(x_1, \nu)]) =0 \} \leq
 \frac{||(b(\cdot,\nu ) - b(\cdot,\alpha)) _+ ||_{L^1}}
{|| b(\cdot,\alpha)  ||_{L^1}} \leq \varepsilon.
\]

The case of death occurring at $\theta _1$ may be
analyzed in the same way.

From inequalities \eqref{vernachlassigen} and 
\eqref{zurechtkommen} we may deduce 

 \begin{equation}
  \sup\limits _{t\in (0;\theta _1] }  dist( \eta (\alpha , t ), \eta (  \alpha _n , t ) )
 \overset{p}{\to} 0, n \to \infty .
  \end{equation}

Proceeding in the same manner we may extend this to 
\begin{equation}
  \sup\limits _{t\in (0;\theta _n] }  dist( \eta (\alpha , t ), \eta ( \alpha _n , t ) )
 \overset{p}{\to} 0, n \to \infty ,
  \end{equation}
particularly because of the strong Markov property of a Poisson
point process. In fact, with high probability
the processes $ \eta (\alpha _n, \cdotp )$ and 
$ \eta (\alpha, \cdotp )$ change
up to time $\theta _n$ in the same way in
the following sense: births occur in the 
same places at the same time moments. Deaths occur at the same time moments,
and when  a
point is deleted from $ \eta (\alpha, \cdotp )$, then its counterpart 
is deleted from $ {\eta } (\alpha _n, \cdotp )$.

Since $\theta _n \to \infty$, we get \eqref{accent}.  $\Box $

\begin{3rmk5} \label{rmrk uniform metric D}
 In fact, we have proved an even stronger statement.
 Namely, take $\alpha _n \to \alpha$. Then there 
 exist processes $(\xi ^{(n)}_t)_{t \in [0;T]}$ such 
 that $(\xi ^{(n)}_t)_{t \in [0;T]} \,{\buildrel d \over =}\,
 (\eta(\alpha _n, t))_{t \in [0;T]}$ and 
 \[
   \sup\limits _{t\in [0;T] }  dist( \eta (\alpha , t ), \xi ^{(n)}_t )
 \overset{p}{\to} 0, \ \ n \to \infty.
 \]
 Thus,
 $Law \{ \eta (\alpha , \cdotp ), \cdotp \in (0;T] \}$
 and $Law \{ \eta (\alpha _n , \cdotp ), \cdotp \in (0;T] \}$
 are close in the space of measures over
 $D_{\CYRG _0} $, even when 
 $D_{\CYRG _0} $ is considered as topological space
 equipped with the \textit{uniform} topology 
(induced by metric $dist$), and not with the Skorokhod 
topology.

\end{3rmk5}

\subsection{The martingale problem}

Now we briefly discuss the martingale problem
associated with $L$ defined in \eqref{the generator}.
Let
$C_b(\Go)$ be the space of all bounded continuous functions on $\Go$.
We equip $C_b(\Go)$ with the supremum norm.

\begin{def mart pr} \label{def mart pr}
 A probability measure $Q$ on $( D_{\CYRG _0} [0;\infty) , \mathscr {B} ( D_{\CYRG _0} [0;\infty) ))$
 is called a solution to the local martingale problem
associated with $L$ if
$$
M^f_t = f(y(t)) - f(y(0)) - \int\limits _0 ^t Lf (y(s-)) ds , \quad \mathscr {I} _t, \ \ 0 \leq t < \infty ,
$$
 is 
 a local martingale for every $f \in C _b(\CYRG _0)$. Here
 $y$ is the coordinate mapping, $y(t)(\omega) = \omega (t)$,
 $\omega \in D_{\CYRG _0} [0;\infty)$,
$\mathscr {I} _t$ is the completion of $\sigma  ( y (s), 0\leq s \leq t )$ under $Q$.

\end{def mart pr}

Thus, we require $M^f$ to be a local martingale under $Q$ with respect
to $\{\mathscr {I} _t\}_{t \geq 0}$. 
Note that $L$ can be considered
as a bounded operator on $C_b(\Go)$.

\begin{3prop1}
 Let $(\eta(\alpha,t) )_{t \geq 0}$ be a solution to \eqref{se}.
 Then for every $f \in C (\CYRG _0)$ 
 the process 
\begin{equation} \label{martingale}
M^f_t = f(\eta(\alpha,t)) - f(\eta(\alpha,t)) - \int\limits _0 ^t Lf (\eta(\alpha,s-)) ds 
\end{equation}
is a local martingale
under $P$ with respect to $\{ \mathscr{F} _t \}_{t \geq 0}$.
\end{3prop1}

\textbf{Proof}. In this proof $\zeta _t$ will stand for $\eta (\alpha, t)$. 
Denote $\tau _n = \inf \{t\geq 0: |\zeta _t| >n 
\text{ or \ } \zeta _t \nsubseteq [-n;n]^d \}$. 
Clearly, $\tau _n$,
$n \in \N$, is a stopping time and $ \tau _n \to \infty $ a.s.
 Let $\zeta ^n _t = \zeta _{t \wedge \tau _n}$. 
We want to show that $(^{(n)}M^f_t)_{t\geq 0}$ is a martingale, where 

\begin{equation}\label{verisimilitude}
^{(n)}M^f_t = f(\zeta ^n _t) - f(\zeta ^n _t) 
- \int\limits _0 ^t Lf (\zeta ^n _{s-}) ds  .
\end{equation}

 The process $(\zeta _t )_{t \geq 0}$
 satisfies

\begin{equation}\label{credence}
\zeta _t  = \sum\limits _{s\leq t, \zeta  _s \ne \zeta _{s-}} 
[\zeta  _s  - \zeta _{s-}] +\zeta _{0}.
\end{equation}
In the above equality as well as in few other 
places throughout this proof we treat elements
of $\Go$ as measures rather than as configurations.
Since $(\zeta _t)$ is of the pure jump type,
the sum on the right-hand side of 
\eqref{credence}
 is a.s. finite.
Consequently 
we have 
\begin{align} \label{mart111}
&f(\zeta ^n _t) - f(\zeta ^n _0) = 
\sum\limits _{s\leq t, \zeta  _s \ne \zeta _{s-}} [f(\zeta ^n _s)  - f(\zeta ^n _{s-})]  \notag \\
= \int\limits _{B \times (0;t] \times [0; \infty ] } &[f(\zeta  _s)  - f(\zeta _{s-})]
I_{\{s\leq \tau _n\}} I _{ [0;b(x,\zeta _{s-} )] } (u) dN_1(x,s,u) \\
- \int\limits _{\Z \times (0;t] \times [0; \infty ] } I _{\{x_i \in \zeta _{s-}\}} &
[f(\zeta _{s})  - f(\zeta _{s-})] I_{\{s\leq \tau _n\}}
I _{ [0;d(x_i,\zeta _{s-} )] } (v) dN_2 (i,s,v) . \notag
\end{align}

Note that $\zeta _s = \zeta _{s-} \cup x$ a.s. in the first 
summand on the right-hand side of \eqref{mart111}, and $\zeta _s = \zeta _{s-} \setminus x_i$
a.s.
in the second summand.
Now, 
we may write
\begin{align} \label{mart222}
&\int\limits _0 ^t I_{\{s\leq \tau _n\}} L f (\zeta  _s) ds = \notag \\
\int\limits _0 ^t \int\limits _{x \in \R^d , u\geq 0} I_{\{s\leq \tau _n\}}&
I _{ [0;b(x,\zeta _{s-} )] } (u) [f(\zeta _{s-}\cup {x}) -f(\zeta _{s-})] dx du ds -  \\
\int\limits _0 ^t  \int\limits _{x \in \R^d , u\geq 0} I_{\{s\leq \tau _n\}}&
 I _{ [0;d(x,\zeta _{s-}) )] } (v) [f(\zeta _{s-} \setminus {x}) - f(\zeta _{s-})] 
 \zeta _{s-}(dx) dv ds. \notag 
\end{align}

Functions $b,d(\cdot , \cdot )$ and $f$ are bounded on 
$\R ^d \times \{ \alpha : |\alpha | \leq n \text{ and } \alpha \subset [-n;n]^d  \} $
and $\{ \alpha : |\alpha | \leq n \text{ and } \alpha \subset [-n;n]^d  \}$
respectively
by a constant $C>0$.
 Now, for a predictable bounded processes $(\gamma _s (x,u)) _{0 \leq s\leq t}$ and
$(\beta _s (x,v)) _{0 \leq s\leq t}$, the processes

\begin{align}
\int\limits _{B \times (0;t] \times [0; C ] } I_{\{s\leq \tau _n\}}\gamma _s (x,u) [&dN_1(x,s,u) - dxdsdu], \notag \\
\int\limits _{\Z \times (0;t] \times [0; C ] } I_{\{s\leq \tau _n\}}
I_{\{x_i \in \zeta _{s-}\}} \beta _s (x_i,v) &[dN_2 (i,s,v) - \#(di) dsdv] . \notag
\end{align}
are martingales. 
Observe that 
\[
\int\limits _{\Z \times (0;t] \times [0; C ] } I_{\{s\leq \tau _n\}}
I_{\{x_i \in \zeta _{s-}\}} \beta _s (x_i,v)  \#(di) dsdv =
\int\limits _{\Z \times (0;t] \times [0; C ] } I_{\{s\leq \tau _n\}}
 \beta _s (x,v)  \zeta _{s-}(dx) dsdv
\]

Taking

\[\gamma _s (x,u) =
I _{ [0;b(x,\zeta _{s-} )] } (u) [f(\zeta _{s-}\cup {x}) -f(\zeta _{s-})],
\]
\[\beta _s (x,v)) =
 I _{ [0;d(x,\zeta _{s-})] } (v) [f(\zeta _{s-} \setminus {x}) - f(\zeta _{s-})],
 \] 
we see that the difference on the right hand side of \eqref{verisimilitude} is a martingale
because of \eqref{mart111} and \eqref{mart222}. $\Box$
 
\begin{3crl mart pr}
The unique solution of \eqref{se} induces a solution of the martingale problem \ref{def mart pr}.
\end{3crl mart pr}

\begin{3rmk7}
 Since $y(s) = y(s-)$ $P_\alpha$ - a.s., the process 
\[
 f(y(t)) - f(y(0)) - \int\limits _0 ^t Lf (y(s)) ds ,  \  0 \leq t < \infty ,
\]
 is a local martingale, too.
\end{3rmk7}

\subsection{Birth rate without sublinear growth condition}

In this section we will consider equation \eqref{se}
with the a birth rate coefficient that does not satisfy
the sublinear growth condition
\eqref{sublinear growth for b}.

Instead, we assume only that

\begin{equation} \label{condition on b}
\sup\limits _{x\in \R^d,|\eta| \leq m} b(x, \eta) < \infty.
\end{equation}

Under this assumption we can not guarantee existence of 
solution on the whole line $[0;\infty)$ or even on a finite
interval $[0;T]$. It is possible that infinitely many points appear in finite time.

We would like to show
that a unique solution exists up to an explosion time, maybe finite. Consider
birth and death coefficients 

\begin{equation}\label{truncated b,d}
\begin{split}
  b_n (x, \eta) = b(x, \eta ) I_{\{ |\eta| \leq n \}}, \\
  d_n (x, \eta) = d(x, \eta ) I_{\{ |\eta| \leq n \}}.
\end{split}
\end{equation}

Functions $b_n, d_n$ are bounded,
so equation \eqref{se} with 
 birth rate coefficient $b_n$
and death rate coefficient 
$d_n$ has a unique solution
by
Theorem \ref{ex un}. 
Remark \ref{3rmk2} provides 
the existence and uniqueness of solution to \eqref{se} (with
 birth and death rate coefficients 
$b$ and $d$, respectively) up
to the (random stopping) time $\tau _n = \inf \{ s\geq 0 : |\eta _s| >n  \}$. 
Clearly, $\tau _{n+1}\geq \tau _n$; 
if $\tau _n \to \infty$ a.s., 
then we have existence and uniqueness for \eqref{se}; if $\tau _n \uparrow \tau < \infty$
with positive probability, then we have an \textit{explosion}. 
However, existence and uniqueness hold 
up to explosion time $\tau$. 
When we have an explosion we say that the solution blows up.

\subsection{Coupling}

Here we discuss the coupling
of two birth-and-death processes. 
The theorem we prove here
will be used in the sequel. As a matter of fact,
we have already used the coupling technique in
the proof of Theorem \ref{ex un}.

Consider two
equations of the form \eqref{se},

\begin{equation} \label{2se}
\begin{split}
\xi ^{(k)} _t (B) &= \int\limits _{B \times (0;t] \times [0; \infty ] }
I _{ [0;b _k(x,\xi ^{(k)} _{s-} )] } (u) dN_1(x,s,u) \\
-  \int\limits _{\Z \times (0;t] \times [0; \infty ) } I&_{\{x^{(k)}_i \in \xi ^{(k)} _{r-} \cap B \}}
I _{ [0;d(x^{(k)}_i,\eta _{r-} )] } (v)  d  N _2 (i,r,v) + \xi ^{(k)} _{0}(B) , \ \  k=1,2, 
\end{split}
\end{equation}
where $t \in [0;T]$ and
$\{ x^{(k)}_i \}$ is the sequence related to $(\xi ^{(k)}_t)_{t \in [0;T]}$.

Assume that initial conditions $\xi ^{(k)}_0$ and coefficients
$b_k$, $d_k$ satisfy the conditions of Theorem \ref{ex un}. Let
$(\xi ^{(k)}_t)_{t \in [0;T]}$ 
be the unique strong solutions.  

\begin{3thm5} \label{couple}
 Assume that almost surely $\xi ^{(1)}_0 \subset \xi ^{(2)}_0$, and
for any two finite configurations $\eta ^1 \subset \eta ^2$,
\begin{equation} \label{culpable}
b_1(x,\eta ^1) \leq b_2 (x, \eta ^2), \ \ \ x\in \R^d
\end{equation}
and 
$$
d_1(x,\eta ^1) \geq d_2 (x, \eta ^2), \ \ \ x\in \eta ^1.
$$

Then there exists a
cadlag $\Go$-valued process $ ( \eta_t)_{t \in [0;T]}$
such that $ ( \eta_t)_{t \in [0;T]}$ and $( \xi ^{(1)}_t)_{t \in [0;T]}$
have the same law and 

\begin{equation} \label{culprit}
 \eta_t \subset \xi ^{(2)}_t, \ \ \ t \in [0;T].
\end{equation}

\end{3thm5}

\textbf{Proof}. Let 
$\{...,x_{-1} ^{(2)},x_0 ^{(2)},x_1 ^{(2)},... \}$ be the sequence related to 
$(\xi ^{(2)}_t)_{t \in [0;T]}$.
Consider the equation
\begin{equation} \label{apprise}
\begin{split}
\eta _t (B) &= \int\limits _{B \times (0;t] \times [0; \infty ] }
I _{ [0;b _k(x,\eta _{s-} )] } (u) dN_1(x,s,u) \\
-  \int\limits _{\Z \times (0;t] \times [0; \infty ) } I&_{\{x^{(2)}_i \in \eta _{r-} \cap B \}}
I _{ [0;d(x^{(2)}_i,\eta _{r-} )] } (v)  d  N _2 (i,r,v) + \xi ^{(1)} _{0}(B) , \ \  k=1,2.
\end{split}
\end{equation}
Note that here $\{ x^{(2)}_i \}$ is related to $(\xi ^{(2)}_t)_{t \in [0;T]}$ 
and not to $(\eta _t)_{t \in [0;T]}$. Thus \eqref{apprise} is
not an equation of form \eqref{se}.
Nonetheless, the existence of a unique solution
can be shown
in the same way as in the proof of Theorem \ref{ex un}.
Denote the unique strong solution
of \eqref{apprise}
by $(\eta_t)_{t \in [0;T]}$.

Denote by $\{ \tau _m \} _{m \in \N}$ the moments of 
 jumps of $(\eta_t)_{t \in [0;T]}$ and $(\xi ^{(2)}_t)_{t \in [0;T]}$,
$0<\tau _1 <\tau _2 < \tau _3 < ...$. More precisely, a time
$t\in \{ \tau _m \} _{m \in \N}$ iff 
at least one of the processes $(\eta_t)_{t \in [0;T]}$ 
and $(\xi ^{(2)}_t)_{t \in [0;T]}$ jumps at time $t$.

 We will show by induction
that each moment of birth for $(\eta_t)_{t \in [0;T]}$ is
a moment of birth for $(\xi ^{(2)}_t)_{t \in [0;T]}$ too, and 
each moment of death for $(\xi ^{(2)}_t)_{t \in [0;T]}$
is a moment of death for $(\eta_t)_{t \in [0;T]}$
if the dying point is in $(\eta_t)_{t \in [0;T]}$.
Moreover, in both cases
the birth or the death occurs at exactly the same point. Here a moment
of birth is a random time at which a new point appears, a
moment of death is a random time at which a point disappears
from the configuration. The statement formulated
above is in fact equivalent
to \eqref{culprit}.

 Here we deal only with
the base case, the induction step is done in the same way.
We have nothing to show 
if $\tau _1$ is a moment of a birth of $(\xi ^{(2)}_t)_{t \in [0;T]}$
or a moment of death of $(\eta_t)_{t \in [0;T]}$.
Assume that a new point is born for  $(\eta_t)_{t \in [0;T]}$
at $\tau _1$,

$$
\eta_{\tau _1} \setminus \eta_{\tau _1 - } =\{x_1 \}.
$$
The process $(\eta_t)_{t \in [0;T]}$ satisfies \eqref{apprise}, therefore
$N_1(\{x \} , \{ \tau _1 \} , [0;b _k(x_1, \eta _{\tau _1-} )]) = 1$.
Since 
$$ 
\eta _{\tau _1-} = \xi ^{(1)} _0 \subset \xi ^{(2)} _0 = \xi ^{(2)} _{\tau _1-},
$$
by \eqref{culpable} 
$$
N_1(\{x \} , \{ \tau _1 \} , [0;b _k(x_1, \xi ^{(2)} _{\tau _1-} )]) = 1,
$$
hence 
$$
\xi ^{(2)}_{\tau _1} \setminus \xi ^{(2)}_{\tau _1 - } =\{x_1 \}.
$$
The case when $\tau _2$ is a moment of 
death for $(\xi ^{(2)}_t)_{t \in [0;T]}$ is analyzed
analogously.

It remains to show that $ ( \eta_t)_{t \in [0;T]}$
and $( \xi ^{(1)}_t)_{t \in [0;T]}$
have the same law. We mentioned above that formally
equation \eqref{apprise} is not of the form \eqref{se},
so we can not directly apply the uniqueness in 
law result. However, since $\eta _t \in  \xi ^{(2)} _t$
a.s., $t \in [0;T]$, we can still consider 
\eqref{apprise} as an equation of the form \eqref{se}.
Indeed, let $\{...,y_{-1},y_0,y_1,... \}$ be the sequence
 related to $\eta _t$. We have 
$\{y_{-|\xi^{(1)}_0|+1},...,y_{-1},y_0,y_1,... \} \subset 
\{x_{-|\xi^{(2)}_0|+1},...,x^{(2)}_{-1},x^{(2)}_0,x^{(2)}_1,... \}$.
There exists an injection
$\varsigma: \{-|\xi^{(1)}_0|+1,...,0,1,... \} \to \{-|\xi^{(2)}_0|+1,...,0,1,... \}$
such that $y_{\varsigma (i)} = x_i$. Denote
$\theta _i =\inf \{s\geq 0 : y_i \in \eta _s\} $.
Note that $\theta _i$ is a stopping time
with respect to $\{\mathscr{F}_t \}$.
Define a Poisson point process $\bar N_2$ by
\[
 \bar N_2(\{i \} \times R \times V ) =  N_2(\{i \} \times R \times V ), \ \ 
 i \in \Z, R \subset [0;\theta _i], V \subset \R _+,
\]
and
\[
 \bar N_2(\{i \} \times R \times V ) =  N_2(\{ \varsigma(i) \} \times R \times V ), \ \ 
 i \in \Z, R \subset (\theta _i; \infty), V \subset \R _+.
\]
The process $\bar N_2$ is $\{\mathscr{F}_t \}$-adapted.
One can see that $ ( \eta_t)_{t \in [0;T]}$
is the unique solution of
equation \eqref{se} with 
$N_2$ replaced by $\bar  N_2$.
Hence $ ( \eta_t)_{t \in [0;T]} \,{\buildrel d \over =}\,
( \xi ^{(1)}_t)_{t \in [0;T]}$.

\subsection{Related semigroup of operators}

We say now a few words about the semigroup
of operators related to the unique solution 
of \eqref{se}.
We write $\eta(\alpha , t)$ for a unique solution of \eqref{se},
started from $\alpha \in \Go$. 
We want to define an operator $S_t$ 
by
\begin{equation}\label{commiserate}
 S_t f (\alpha) = E f(\eta(\alpha , t)) \ \ \ (= E _{\alpha} f(\eta( t)))
\end{equation}
for 
an appropriate class of functions.
Unfortunately, it seems difficult to make $S_t$
a $C_0$-semigroup on some functional Banach space
for general $b,d$ satisfying \eqref{sublinear growth for b}
and \eqref{condition on d}.

We start with the case 
when the cumulative birth and death 
rates are bounded. Let $C_b=C _b (\Go)$
be the space of all bounded continuous functions on $\Go$.
It becomes a Banach space once it is 
equipped with the supremum norm.
We assume the existence of a constant $C>0$
such that
for all 
$\zeta \in \Go$ 
\begin{equation}\label{dabble at}
 |B(\zeta)|+|D(\zeta)| <C,
\end{equation}
where
$B$ and $D$ are defined in \eqref{cumulative death rate} and
\eqref{cumulative birth rate}.
Formula \eqref{the generator} defines
then a bounded operator $L: C_b \to C_b$, 
and we will show that $S_t$ coincides with $e^{tL}$.
For $f\in C_b$, the function $S_t f$
is bounded and continuous.
Boundedness is a consequence
 of the boundedness of $f$, and  continuity of $S_t f$
 follows from Remark \ref{rmrk uniform metric D}.
 Indeed, let $\alpha _n \to \alpha$, $\xi ^{(n)} _t \overset{d}{=} \eta (\alpha _n ,t)$
 and 
  \[
    dist( \eta (\alpha , t ), \xi ^{(n)}_t )
 \overset{p}{\to} 0, \ \ n \to \infty.
 \]
 Unlike $\G$, the space $\Go$ is a $\sigma$-compact space. 
 Consequently, for all $\varepsilon >0$ there exists a compact
 $K_{\varepsilon} \subset \Go$ such that for large enough $n$ 
 \[
  P \{ \eta (\alpha , t ) \in K_{\varepsilon} , \  \xi ^{(n)} _t \in K_{\varepsilon}
  \} \geq 1 - \varepsilon.
 \]
 Also, for fixed $\delta >0$ and for large enough $n$
 \[
  P \{ 
  dist( \eta (\alpha , t ), \xi ^{(n)}_t) \leq \delta \} \geq 1 - \delta.
 \]
 
Fix $\varepsilon >0$. There exists $\delta _{\varepsilon} \in (0;\varepsilon)$ 
such that $|f(\beta) - f (\gamma)| \leq \varepsilon $
whenever $dist( \beta, \gamma) \leq \delta _{\varepsilon}$,
$\beta, \gamma \in K_{\varepsilon}$. 
We have for large enough $n$
\[
 |E[f(\eta (\alpha , t ) ) - f(\xi ^{(n)}_t)]|
 \]
 \[\leq 
 E|f(\eta (\alpha , t ) ) - f(\xi ^{(n)}_t)|
 I\{\eta (\alpha , t ) \in K_{\varepsilon} , \  \xi ^{(n)} _t \in K_{\varepsilon},
 dist( \eta (\alpha , t ), \xi ^{(n)}_t) \leq \delta _{\varepsilon} \} 
 \]
 \[+ 
 2 (\delta _{\varepsilon} + \varepsilon) ||f|| \leq
 \varepsilon + 2 (\delta _{\varepsilon} + \varepsilon) ||f||,
\]
where $||f|| = \sup _{\zeta \in \Go} |f(\zeta)|$.
Letting $\varepsilon \to 0$, we see that 
\[
 E f(\eta (\alpha _n , t ) ) = E f(\xi ^{(n)}_t) \to E f(\eta (\alpha, t ) ).
\]
 Thus, $S_t f$ is continuous
 (note that the continuity of $S_t f$  does not follow
 from Theorem \ref{continuity of solution}
 alone, since for a fixed $t\in [0;T]$
 the functional 
 $D_{\Go} [0;T] \ni x \mapsto x(t) \in \R$
 is not continuous in the  Skorokhod topology).
 Furthermore, since for small $t$ and for all ${A \in \mathscr{B}(\R ^d)}$,
 
 \begin{equation}\label{verheimlichen}
  P \{ \eta (\alpha, t) = \alpha \} = 1 - 
  t [B(\alpha) + D(\alpha)] + o(t),
 \end{equation}
 
 \begin{equation}\label{redound}
  P \{ \eta (\alpha, t) = \alpha \cup \{ y \} \text{ for some }
  y \in A  \} = t \int\limits _{y \in A} b(y, \alpha) dy
  + o(t),
 \end{equation}
and for $x \in \alpha$
\begin{equation}\label{sanctity}
  P \{ \eta (\alpha, t) = \alpha \setminus \{x\}  \} = t d(x,\alpha) 
  + o(t),
 \end{equation}

 we may estimate
 
 \[
  |S_t f ( \alpha ) - f(\alpha)| \leq t \left[B(\alpha) + D(\alpha) \right] ||f||
  +o(t) ||f||\leq C ||f|| t +o(t).
 \]
Therefore, \eqref{commiserate} defines 
a $C_0$ semigroup on $C_b$.
Its generator 

\[
 \tilde L f(\alpha) = \lim _{t \to 0+}\frac{S_t f(\alpha)}{t} = 
\]
\[
\lim\limits _{t \to 0+} \left[ \int\limits _{x \in \R^d} b(x, \alpha) [f(\alpha \cup {x}) -f(\alpha)] dx + 
\sum\limits _{x \in \alpha} d(x, \alpha ) (f(\alpha \setminus {x}) - f(\alpha))  +o(t)\right]
=L f(\alpha).
\]

Thus, $ S_t = e^{tL}$, and we have proved the following

\begin{8prop1}
 Assume that \eqref{dabble at} is fulfilled. Then the family
 of operators $(S_t, t \geq 0)$ on $C_b$ defined in \eqref{commiserate}
 constitutes a $C_0$-semigroup. Its generator coincides with $L$ given in 
 \eqref{the generator}.
 
\end{8prop1}

Now we turn out attention to general $b,d$
satisfying \eqref{sublinear growth for b}
and \eqref{condition on d} but 
not necessarily \eqref{dabble at}.
The family of operators $(S_t)_{t\geq 0}$
still constitutes a semigroup,
however it does not 
have to be strongly continuous 
anymore.
Consider truncated birth
and death coefficients \eqref{truncated b,d}
and corresponding process $\eta ^n (\alpha, t)$.
Remark \ref{3rmk2} implies
that $\eta ^n (\alpha, t)= \eta  (\alpha, t)$
for all $t \in [0;\tau_n]$, where

\begin{equation}
 \tau _n = \inf \{ s \geq 0 : | \eta (\alpha,s)| > n \}.
\end{equation}

Growth condition \eqref{sublinear growth for b} implies that
$\tau _n \to \infty$ for any $\alpha \in \Go$.

Truncated coefficients $b_n, d_n$ satisfy 
\eqref{dabble at} and

\begin{equation}
 S_t ^{(n)} f (\alpha) = E f(\eta ^{(n)}(\alpha , t))
\end{equation}
defines
a $C_0$ - semigroup on $C_b$. In particular, 
for all $\alpha \in \Go$
\[
 L ^{(n)} f (\alpha) = \lim\limits _{t \to 0+} 
 \frac{E f(\eta ^{(n)}(\alpha , t)) - f (\alpha)}{t},
\]
where $L ^{(n)}$ is 
operator defined as in \eqref{the generator}
but with $b_n,d_n$ instead of $b,d$.
Letting $n \to \infty$ we get, for fixed $\alpha$
and $f$, 

\begin{equation}\label{gloss}
 L  f (\alpha) = \lim\limits _{t \to 0+} 
 \frac{E f(\eta (\alpha , t)) - f (\alpha)}{t} = \lim\limits _{t \to 0+} 
 \frac{S_t f(\alpha) - f (\alpha)}{t}.
\end{equation}

Taking limit by $n$ is possible: for $n \geq |\alpha| +2$,
 $\eta ^{(n)}(\alpha , t)$
satisfies \eqref{verheimlichen}, \eqref{redound} and 
\eqref{sanctity}, 
therefore $\eta (\alpha , t)$ satisfies \eqref{verheimlichen}, \eqref{redound} and 
\eqref{sanctity}, too.
Thus, we have 

\begin{8prop2}
 Let $b$ and $d$ satisfy \eqref{sublinear growth for b} and 
 \eqref{condition on d} but not necessarily \eqref{dabble at}.
 Then the family of operators $(S_t, t \geq 0)$ constitutes
 a semigroup on $C_b$ which does not have to be strongly continuous.
 However,
 for every $\alpha \in \Go$ and $f \in C_b$ we have \eqref{gloss}.
 
\end{8prop2}

Formula \eqref{gloss} gives us the formal relation of 
 $(\eta (\alpha , t))_{t \geq 0}$ to the operator $L$.
Of course, for fixed $f$ the convergence in \eqref{gloss}
does not have to be uniform in $\alpha$.

\begin{8rmk10}
 The question about the construction of a semigroup 
 acting on some class of probability
 measures on $\Go$ is 
 yet to be studied.

\end{8rmk10}

  \section*{Acknowledgements}

 This paper is part of the author's 
 PhD thesis written under the supervision
 of Professor Yuri Kondratiev,
 whom the author would like to thank for suggesting the problem
 and stimulating discussions.
  The financial support of the German science foundation
  through the IGK is gratefully appreciated.

\bibliographystyle{alpha}
\bibliography{BaDdynamics}

\newcommand{\etalchar}[1]{$^{#1}$}
\begin{thebibliography}{FOK{\etalchar{+}}14}

\bibitem[Ald13]{IPS1}
D.~Aldous.
\newblock Interacting particle systems as stochastic social dynamics.
\newblock {\em Bernoulli}, 19:1122--1149, 2013.

\bibitem[AN72]{Branch1}
K.B. Athreya and P.E Ney.
\newblock {\em Branching processes}.
\newblock Die Grundlehren der Mathematischen Wissenschaften in
  Einzeldarstellungen. Springer, 1972.

\bibitem[Arn06]{Yuleproc}
d.L.F. Arnaud.
\newblock Yule process sample path asymptotics.
\newblock {\em Electron. Comm. Probab.}, 11:193–199, 2006.

\bibitem[BR58]{BRosenblatt}
C.~J. Burke and M.~Rosenblatt.
\newblock A {M}arkovian function of a markov chain.
\newblock {\em E Ann. Math. Statist}, 29:1112–1122, 1958.

\bibitem[CG]{Onaseminal}
M.~M. Castro and F.~A. Grünbaum.
\newblock On a seminal paper by karlin and mcgregor.
\newblock {\em Symmetry Integrability Geom. Methods Appl.}, 9.

\bibitem[EK86]{EthierKurtz}
S.~N. Ethier and T.~G. Kurtz.
\newblock {\em Markov Processes. Characterization and convergence}.
\newblock Wiley-Interscience, New Jersey, 1986.

\bibitem[FKK12a]{BaDdynamics}
D.~Finkelshtein, O.~Kutovyi, and Yu. Kondratiev.
\newblock Semigroup approach to birth-and-death stochastic dynamics in
  continuum.
\newblock {\em Journal of Functional Analysis}, 262(3):1274–1308, 2012.

\bibitem[FKK12b]{Semigroupapproach}
D.~Finkilstein, Yu. Kondratiev, and O.~Kutoviy.
\newblock Semigroup approach to birth-and-death stochastic dynamics in
  continuum.
\newblock {\em J. Funct. Anal.}, 262(3):1274–1308, 2012.

\bibitem[FKK14]{Statdynamics}
D.~Finkelshtein, O.~Kutovyi, and Yu. Kondratiev.
\newblock Statistical dynamics of continuous systems: perturbative and
  approximative approaches.
\newblock {\em Arabian Journal of Mathematics}, 2014.
\newblock doi:10.1007/s40065-014-0111-8.

\bibitem[FM04]{FournierMeleard}
N.~Fournier and S.~Méléard.
\newblock A microscopic probabilistic description of a locally regulated
  population and macroscopic approximations.
\newblock {\em Ann. Appl. Probab}, 14(4):1880–1919, 2004.

\bibitem[FOK{\etalchar{+}}14]{Ecology}
D.~Finkelshtein, O.~Ovaskainen, O.~Kutovyi, S.~Cornell, B.~Bolker, and Yu.
  Kondratiev.
\newblock A mathematical framework for the analysis of spatial-temporal point
  processes.
\newblock {\em Theoretical Ecology}, 7:101--113, 2014.

\bibitem[Fra14]{IPS2}
T.~Franco.
\newblock Interacting particle systems: hydrodynamic limit versus high density
  limit.
\newblock 2014.
\newblock preprint; arXiv:1401.3622 [math.PR].

\bibitem[Gar95]{Garcia}
N.~L. Garcia.
\newblock Birth and death processes as projections of higher-dimensional
  poisson processes.
\newblock {\em Adv. in Appl. Probab.}, 27(4):911–930, 1995.

\bibitem[GK06]{GarciaKurtz}
N.~L. Garcia and T.~G. Kurtz.
\newblock Spatial birth and death processes as solutions of stochastic
  equations.
\newblock {\em ALEA Lat. Am. J. Probab. Math. Stat.}, (1):281–303, 2006.

\bibitem[GK08]{GarciaKurtz2}
N.~L. Garcia and T.~G. Kurtz.
\newblock Spatial point processes and the projection method.
\newblock {\em Progr. Probab. In and out of equilibrium. 2,}, 60(2):271–298,
  2008.

\bibitem[GS75]{GikhSkor2}
I.~I. Gikhman and A.~V. Skorokhod.
\newblock {\em The Theory of Stochastic Processes}, volume~2.
\newblock Springer, 1975.

\bibitem[GS79]{GikhSkor3}
I.~I. Gikhman and A.~V. Skorokhod.
\newblock {\em The Theory of Stochastic Processes}, volume~3.
\newblock Springer, 1979.

\bibitem[Har63]{Branch2}
T.E. Harris.
\newblock {\em The theory of branching processes}.
\newblock Die Grundlehren der Mathematischen Wissenschaften in
  Einzeldarstellungen. Springer, 1963.

\bibitem[HS78]{HolleyStroock}
R.~A. Holley and D.~W. Stroock.
\newblock Nearest neighbor birth and death processes on the real line.
\newblock {\em Acta Math}, 140(1-2):103–154, 1978.

\bibitem[IW81]{IkedaWat}
N.~Ikeda and S.~Watanabe.
\newblock {\em Stochastic Differential Equations and Diffusion Processes}.
\newblock Nord-Holland publiching company, 1981.

\bibitem[Kal02]{KallenbergFound}
O.~Kallenberg.
\newblock {\em Foundations of modern probability}.
\newblock Springer, 2 edition, 2002.

\bibitem[Kin93]{KingmanPP}
J.~F.~C. Kingman.
\newblock {\em Poisson Processes}.
\newblock Oxford University Press, 1993.

\bibitem[KK02]{KondKuna}
Y.~G. Kondratiev and T.~Kuna.
\newblock Harmonic analysis on configuration space. i. general theory.
\newblock {\em Infin. Dimens. Anal. Quantum Probab. Relat. Top.},
  5(2):201--233, 2002.

\bibitem[KK06]{KondKut}
Yu.~G. Kondratiev and O.~V. Kutoviy.
\newblock On the metrical properties of the configuration space.
\newblock {\em Math. Nachr.}, 279:774–783, 2006.

\bibitem[KL99]{Scalinglimits}
C.~Kipnis and C.~Landim.
\newblock {\em Scaling limits of interacting particle systems}.
\newblock Springer, 1999.

\bibitem[KM59]{KarlinMcGregor}
S.~Karlin and J.~McGregor.
\newblock Random walks.
\newblock {\em Illinois J. Math.}, 3:66--81, 1959.

\bibitem[KS06]{KondSkor}
Yu. Kondratiev and A.~Skorokhod.
\newblock On contact processes in continuum.
\newblock {\em Infin. Dimens. Anal. Quantum Probab. Relat. Top.},
  9(2):187–198, 2006.

\bibitem[Kur66]{Kuratowski}
K.~Kuratowski.
\newblock {\em Topology}, volume~1.
\newblock Academic Press, New York and London, 1966.

\bibitem[Lev03]{Levin}
S.A. Levin.
\newblock Complex adaptive systems:exploring the known, the unknown and the
  unknowable.
\newblock {\em Bulletin of the AMS}, 40(1):3--19, 2003.

\bibitem[Lig85]{Liggett}
T.~M. Liggett.
\newblock {\em Interacting particle systems}.
\newblock Grundlehren der Mathematischen Wissenschaften. Springer, 1985.

\bibitem[Lig04]{Liggett2}
T.~M. Liggett.
\newblock {\em Interacting particle systems—an introduction}.
\newblock 2004.
\newblock ICTP Lect. Notes, XVII.

\bibitem[MS94]{Statsim}
J.~Møller and M.~Sørensen.
\newblock Statistical analysis of a spatial birth-and-death process model with
  a view to modelling linear dune fields.
\newblock {\em Scand. J. Statist.}, 21(1):1--19, 1994.

\bibitem[MW04]{Simbook}
J.~Møller and R.~P. Waagepetersen.
\newblock {\em Statistical Inference and Simulation for Spatial Point
  Processes}.
\newblock Chapman and Hall/CRC, 2004.

\bibitem[Pen08]{Penrose}
M.~D. Penrose.
\newblock Existence and spatial limit theorems for lattice and continuum
  particle systems.
\newblock {\em Probab. Surv.}, 5:1--36, 2008.

\bibitem[Pod09]{Uncountable}
K.~Podczeck.
\newblock On existence of rich fubini extensions.
\newblock {\em Econom. Theory}, 45(1-2):1--22, 2009.

\bibitem[Pre75]{Preston}
C.~Preston.
\newblock Spatial birth-and-death processes.
\newblock In {\em Proceedings of the 40th Session of the International
  Statistical Institute}, volume~46 of {\em Bull. Inst. Internat. Statist},
  pages 371–391, 405–408, 1975.

\bibitem[RS99]{RockSchied}
M.~Röckner and A.~Schied.
\newblock Rademacher's theorem on configuration spaces and applications.
\newblock {\em J. Funct. Anal.}, 169(2):325–356, 1999.

\bibitem[RY05]{RevuzYor}
D.~Revuz and M.~Yor.
\newblock {\em Continuous Martingales and Brownian Motion}.
\newblock Springer, 3 edition, 2005.

\bibitem[Spi77]{SpitzerBaD}
F.~Spitzer.
\newblock Stochastic time evolution of one dimensional infinite particle
  systems.
\newblock {\em Bull. Amer. Math. Soc.}, 83(5):880–890, 1977.

\end{thebibliography}

\end{document}